\documentclass[11pt]{amsart}
\usepackage{palatino}
\usepackage[all]{xy,xypic}
\usepackage{amsfonts,amsmath,oldgerm,amssymb,amscd,tikz-cd,enumerate}
\usepackage [mathscr]{euscript}
\makeatletter
\@namedef{subjclassname@2020}{%
\textup{2020} Mathematics Subject Classification}
\makeatother

\newcommand{\ra}{\rightarrow}		
\newcommand{\lra}{\longrightarrow}
\newcommand{\by}[1]{\stackrel{#1}{\ra}}

\newcommand{\remove}[1]{}

\newcommand{\surj}{\ra\!\!\!\ra}	

\newcommand{\ol}{\overline}		
\newcommand{\wt}{\widetilde}
\newcommand{\iso}{\by \sim}

\newtheorem{theorem}{Theorem}[section]
\newtheorem{proposition}[theorem]{Proposition}
\newtheorem{lemma}[theorem]{Lemma}
\newtheorem{definition}[theorem]{Definition}
\newtheorem{corollary}[theorem]{Corollary}
\newtheorem{conjecture}[theorem]{Conjecture}
\newtheorem{example}[theorem]{Example}

\newcommand{\ga}{\alpha}	\newcommand{\gb}{\beta}

\newcommand{\BC}{\mbox{$\mathbb C$}}

\newcommand{\bq}{\mbox{$\mathbb Q$}}	\newcommand{\br}{\mbox{$\mathbb R$}}

\newcommand{\MI}{\mbox{$\mathfrak I$}}

\newcommand{\E}{\text{E}}
\newcommand{\ot}{\mbox{\,$\otimes$\,}}	
\newcommand{\op}{\mbox{$\oplus$}}
\newcommand{\Spec}{\text{Spec}}	
\newcommand{\hh}{\text{ht}}
\newcommand{\rank}{\text{rank}}
\newcommand{\sr}{\text{sr}}

\newcommand{\diag}{\mbox{\rm diag\,}}	
\newcommand{\Aut}{\mbox{\rm Aut\,}}	
\newcommand{\Hom}{\text{Hom}\,}	
\newcommand{\dd}{\text{dim}}	
	
\newcommand{\Um}{\text{Um}}		
\newcommand{\SL}{\text{SL}}
\newcommand{\GL}{\text{GL}}

\newcommand{\sur}{\twoheadrightarrow}
\newcommand{\bp}{\begin{proposition}}
\newcommand{\ep}{\end{proposition}}
\newcommand{\bl}{\begin{lemma}}
\newcommand{\el}{\end{lemma}}
\newcommand{\bt}{\begin{theorem}}
\newcommand{\et}{\end{theorem}}
\newcommand{\bc}{\begin{corollary}}
\newcommand{\ec}{\end{corollary}}
\newcommand{\bd}{\begin{definition}}
\newcommand{\ed}{\end{definition}}
\newcommand{\bco}{\begin{conjecture}}
\newcommand{\eco}{\end{conjecture}}
\newcommand{\bma}{\begin{bmatrix}}
\newcommand{\ema}{\end{bmatrix}}

\renewcommand{\k}{\overline{\mathbb F}_p}

\def\rmk{\refstepcounter{theorem}\paragraph{{\bf Remark} \thetheorem}}
\def\proof{\paragraph{Proof}}

\def\quest{\refstepcounter{theorem}\paragraph{{\bf Question} \thetheorem}}

\title [On a question of Nori]{On a question of Nori:\\ obstructions, improvements, and applications }
\author{Sourjya Banerjee and Mrinal Kanti Das}
\address{(Sourjya Banerjee) Stat-Math Unit, Indian Statistical Institute, 203, B.T. Road, Kolkata - 700108}
\curraddr{ Department of Mathematics and Statistics, 
	Indian Institute of Science Education and Research Kolkata, Campus Road, Mohanpur, West Bengal - 741246}
\email{sourjya.pdf@iiserkol.ac.in, sourjya91@gmail.com}

\address{(Mrinal Kanti Das) Stat-Math Unit, Indian Statistical Institute, 203, B.T. Road, Kolkata - 700108} 
\email{mrinal@isical.ac.in, das.doublelife@gmail.com}

\keywords{Efficient generation of ideals, complete intersection, projective modules}
\subjclass[2020]{13C10, 19A15}

\begin{document}
\begin{abstract}
This article concerns a question asked by M. V. Nori on homotopy of sections of projective modules defined on the polynomial algebra over a smooth affine domain $R$. While this question has an affirmative answer, it is known that the assertion does not hold if: (1)
$\dd(R)=2$; or (2) $d\geq 3$ but $R$ is not smooth. We first prove that an affirmative answer can be given for $\dd(R)=2$ when $R$ is an $\k$-algebra. Next, for $d\geq 3$ we find the precise obstruction for the failure in the singular case. 
Further, we improve a result of Mandal (related to Nori's question) in the case when the ring $A$ is an affine $\k$-algebra of dimension $d$. We apply this improvement to 
 define the $n$-th Euler class group $E^n(A)$, where $2n\ge d+2.$ Moreover, if $A$ is smooth, we associate to a unimodular row $v$  of length $n+1$ its Euler class  $e(v)\in E^n(A)$ and show that the corresponding stably free module, say, $P(v)$ has a unimodular element if and only if $e(v)$ vanishes in $E^n(A)$. 
\end{abstract}

\maketitle

\section{Introduction}
Let $X=\Spec(A)$, be a smooth affine variety of dimension $d$. Let $P$ be a projective $A$-module of rank $n$ and $\phi_0:P\twoheadrightarrow I_0$ be a surjective homomorphism. Assume that the zero set of $I_0$, $V(I_0)=Y$ is a smooth affine subvariety of $X$ of dimension $d-n$ and $Z=V(I)$ is a smooth closed subvariety of $X\times \mathbb A^1=\Spec(A[T])$ such that $Z$ intersects $X\times \{0\}$ transversally in $Y\times \{0\}$. In this set up Nori asked the following question (for motivation see \cite[Appendix]{SM}). \\

\quest
Does there exist a surjective map $\phi:P[T]\twoheadrightarrow I/I^2$, which is   compatible with $\phi_0$, have a surjective lift  $\psi:P[T]\twoheadrightarrow I$ such that \\
(i) $\psi|_{T=0}=\phi_0$ and \\
(ii) $\psi|_Z=\phi$ ?\\

In simple commutative algebraic terms, the question of Nori can be restated as follows (we take $P$ to be free for simplicity):\\

\quest\label{Norifree}
Let $A$ be a smooth affine domain of dimension $d$ over an infinite (perfect) field $k$. Let $I\subset A[T]$ be an ideal of height $n$ such that $I=(f_1,\cdots,f_n)+(I^2T)$, where $2n\geq d+3$. Do there exist $g_1,\cdots,g_n\in I$ such that $I=(g_1,\cdots,g_n)$ with $g_i-f_i\in (I^2T)$?\\

This question has been answered affirmatively by S. Mandal \cite[Theorem 2.1]{SM}, when $I$ contains a monic polynomial even without any smoothness assumptions. The local version of the above question was answered positively  by Mandal-Varma \cite{mv}. It was finally settled in the affirmative by Bhatwadekar-Keshari \cite{BHMK}. Prior to that, Bhatwadekar-Sridharan gave an affirmative answer when $n=d\geq 3$. In that same paper of Bhatwadekar-Sridharan, there are two examples which show that Nori's question does not have an affirmative answer if: (1) $n=d=2$; (2) $A$ is not smooth.
 
In Sections 3 and 4 of this paper we deeply investigate the examples given in  \cite{BR} and slice out the reasons for such anomalies. For instance, in case of dimension two, their example of a $\mathbb C$-algebra totally depends on the fact that $SK_1$ of a curve over $\mathbb C$ may be non-trivial. This is in complete contrast to curves over $\k$ (the group $SK_1$ is trivial in this case) and we answer Nori's question in the affirmative when $k=\k$ and $d=n=2$ with certain smoothness assumptions. 

On the other hand, for singular $A$, if $I$ does not contain a monic polynomial, we show that the precise obstruction to Nori's vision distills down to the fact that whether $I\cap A$ is contained only in smooth maximal ideals or not. Please see Section 4 for further details. We also derive various applications of our result in this section which we believe, contribute significantly to the understanding of the behaviour of projective modules and locally complete intersection ideals. For example, we prove (Theorem 
\ref{unimod}):

\bt
Let $R$ be  an affine domain of dimension $d\geq 3$ over an infinite  field $k$ and $P$ be a projective $R[T]$-module of rank $d$. Assume that $P/TP$ has a unimodular element and that there is a surjection $\alpha:P\sur I$, where $I\subset R[T]$ is an  ideal of height $d$ such that $I\cap R$ is contained only in smooth maximal ideals. Then $P$ has a unimodular element.
\et 

Further, with the same $R$ as in the above theorem,  we consider an ideal $I\subset R[T]$ of height $d$ such that $I/I^2$ is generated by $d$ elements and $J:=I\cap R$ is contained only in smooth maximal ideals. In  Theorem \ref{projgen} we roughly  prove that $I$ is surjective image of a projective $R[T]$-module of rank $d$ if and only if $I(0)$ is surjective image of a  projective $R$-module of rank $d$.

As we mentioned at the beginning, Mandal answered Question \ref{Norifree} when $I$ contains a monic polynomial and $2n\geq d+3$. In Section 5 we improve Mandal's theorem when the base field is $\k$. Namely, we prove an analogous result when $2n\geq d+2$ (this includes the case $d=n=2$ which does not hold in general). As a natural extension, we also improve the "relative" version of Mandal's theorem proved by Mandal-Sridharan (see Section 6). These results are crucial to the development of the obstruction theory for stably free $R$-modules of rank $=n\geq \frac{d+2}{2}$ to have a free summand of rank one or not. For an affine $\k$-algebra $R$ of dimension $d$, we first define the $n$-th Euler class group $E^n(R)$. Let $P$ be a stably free module of rank $n$, which is induced by a unimodular 
row $v$ of length $n+1$ over $R$. We associate its Euler class $e(v)\in E^n(R)$ and prove that $e(v)=0$ if and only if $P\simeq Q\oplus R$ for some $R$-module $Q$.  We carry all these out in Sections 7 and 8.

One of the fascinating questions in this area is the so called "monic inversion principle", which roughly entails to asking whether a statement on $A[T]$, if found true after inverting a monic, is true over $A[T]$ to begin with. In this context, we recall a question of Roitman:

\quest
Let $R$ be a commutative Noetherian ring and $P$ be a projective $R[T]$-module. Assume that $P_f$ has a unimodular element, where $f\in R[T]$ is a monic polynomial. Does $P$ have a unimodular element?

Bhatwadekar-Sridharan \cite{BRS01} answered this in the affirmative when $R$ contains an infinite field $k$ and $\dd(R)=\rank (P)=d$. However, much progress has not been made for the cases when $\rank (P) < \dd (R)$. Based on our methods developed in the previous sections, we prove the following result.

\bt
Let $R$ be an affine algebra of dimension $d\geq 2$ over $\k$ and $P$ be a $1$-stably free $R[T]$-module of rank $d-1$. Assume that $P_f$ has a unimodular element, where $f\in R[T]$ is a monic polynomial. Then $P$ also has a unimodular element.
\et

\section{Preliminaries}

All modules considered in this article are assumed to be finitely generated. The purpose of this section is to give a detailed proof of a slightly improved version of a result by Mandal-Murthy \cite[Theorem 3.2]{MM}.  
The next lemma allows us to reduce the proof of Theorem \ref{IMMTT} to the case when the ring is reduced. The proof is essentially contained in \cite[Theorem 2.4]{MKD} hence we opt to skip the proof.
\bl \label{MQ}
Let $A$ be a Noetherian ring, $P$ be a projective $A$-module, $I\subset A$ and $K\subset I^2$ be two ideals. Moreover assume that $\tilde{\phi}:P\surj I/K$ is a surjective map. Suppose that there exists a surjective map $\phi:P\twoheadrightarrow I/I\cap\mathfrak{n}$, which satisfies $\phi\otimes A/ K\cap \mathfrak{n}=\tilde{\phi}\otimes A/ \mathfrak{n}$, where  $\mathfrak{n}$ is the nil-radical of $A$. Then there exists a surjective map $\Phi :P\surj I$ such that $\Phi\otimes A/K=\tilde{\phi}.$
\el

Before going to the next result we recall the following definitions.
\bd Let $A$ be a ring with $1$ and $P$ be a projective $A$-module.
\begin{enumerate}

		\item A row vector $(a_0,...,a_n)\in A^{n+1}$ of length $n+1$ is said to be a unimodular row of length $n+1$, if there exists $(b_0,...,b_n)\in A^{n+1}$ such that $\sum_{i=0}^{n}a_ib_i=1.$ We will denote $\Um_{n+1}(A)$ as the set of all unimodular row vectors of length $n+1$ over the ring $A$.
		\item The stable range (denoted by $\sr(A)$) of $A$ is the smallest natural number $r$ with the property that for any $(u_1,...,u_{r+1})\in \Um_{r+1}(A)$, there exist $\lambda_i\in A$ such that $(u_1+\lambda_1u_{r+1},...,u_r+\lambda_ru_{r+1})\in \Um_r(A)$.
		\item $P$ is said to have a unimodular element if there exists an $A-$linear surjection $\phi:P\surj A$. In this case an element $v\in P$ is said to be a unimodular element of $P$ if $v\in \phi^{-1}(1)$. The set of all unimodular elements of $P$ is denoted by $\Um(P)$.
	\item We say that the projective stable range of $A$ (notation: $psr(A)$) is $n$ if $n$ is the least positive integer such that for any projective $A-$module $P$ of rank $n$ and $(p,a)\in \Um(P\oplus A)$, there exists $q\in P$ such that $p+aq\in \Um(P)$.
	\item Let $P$ be such that $\Um(P)\not=\phi$. We choose $\phi\in P^*$ and $p \in P$ such that $\phi(p)=0$. We define an endomorphism $\phi_p$ as the composite $\phi_p:P\to A\to P$, where $A\to P$ is the map sending $1\to p.$ Then by a transvection we mean an automorphism of $P$, of the form $1+\phi_p$, where either $\phi\in \Um(P^*)$ or $p\in \Um(P)$. By $\E(P)$  we denote  the subgroup of $\Aut(P)$ generated by all transvections.
\end{enumerate}
\ed
\smallskip

The next result can be found in \cite[Theorem 3.7]{AMM}. Here we just restate their result with a slight improvement in the dimension two case. For the proof we just mimic their arguments.
\bl\label{psr}
Let $R$ be an affine algebra of dimension $d\geq 2$ over $\k$ and $\mathfrak{a} \subset R$
be an ideal. Suppose that $P$ is a projective $R$-module of rank $d$ having a unimodular element and $p\in P$ is such that 
$\ol{p}\in \Um(P/\mathfrak{a}P)$, where `bar' denote going modulo $\mathfrak{a}$. Then there exists $q\in \Um(P)$ such that $p\equiv q$ (modulo $\mathfrak{a}$).
\el

\proof
Since transvections have lift, it is enough to have $psr(R)\leq d$. We elaborate.  Since $\ol{p}\in \Um(P/\mathfrak{a}P)$, there exists $a\in \mathfrak{a}$ such that $(p,a)\in \Um(P\op R)$. Now if
$psr(R)\leq d$, then there exists $y\in P$ such that $p+ay\in \Um(P)$. We can then take $q=p+ay$. 

Case $d\ge 3$ follows from \cite[Theorem 3.7]{AMM}. For $d=2$, same proof works if we show that rank $2$ projective $R-$module has unimodular element, which follows from \cite[Theorem 2.5 and Corollary 2.7]{MKDIMRN}. \qed

\medskip

We shall end this section with the following improvement of \cite[Theorem 3.2]{MM}. We essentially follow their proof with some small adjustments to suit our requirements. Before that we shall recall the following definition.

\bd (Order ideal) Let $A$ be a ring. Let $P$ be a projective module and $m\in P$. Then the \textit{order ideal} of $m$ is defined as:
$$\mathcal O(m,P)=\mathcal O(m)=\{f(m):f\in \Hom_A(P,A)\}$$ 
\ed

The next theorem improves upon \cite[Theorem 3.2]{MM} to meet our requirements. The proof is an exact replication of their arguments, so we omit it. This has also been demonstrated by Keshari-Tikader \cite[Theorem 4.1]{STMKK}.

\bt \label{IMMTT}
Let $R$ be an affine algebra of dimension $d\geq 2$ over $\k$ and $I\subset R$ be an ideal. Let $P$ be a projective $R$-module of rank $\ge d$ such that $\ol{f}:P\sur I/K$ be a surjection, where $K\subseteq I^2$ is an ideal. Then $\ol{f}$ lifts to a surjection $f:P\sur I$.
\et

We end this section with the following basic result.

\bp\label{limit}
Let $k$ be a field and $A$ be an affine $k$-algebra of dimension $d$. Let $S$ be a multiplicative subset of $A$ such that $\dd(S^{-1}A)=d$. Then $A$ is the direct limit of affine $k$-algebras, each of dimension $d$.
\ep

\proof
We first note that $S$ is a directed set where, for $s_1,s_2\in S$, we say $s_1\leq s_2$ if and only if there is $s_3\in S$ such that
$s_2=s_1s_3$. It then follows that $\{A_s\,|\,s\in S\}$ is a direct system of rings, indexed by $S$. For $s\in S$, note that 
$$A_s=\frac{A[X]}{\langle sX-1\rangle}.$$
Therefore, each $A_s$ is an affine $k$-algebra of dimension $d$. It is also easy to see that $S^{-1}A=\varinjlim \{A_s\,|\,s\in S\}$.
\qed

\section{Addressing the failures: Dimension two}
In \cite[Example 3.15]{BR}, Bhatwadekar-Sridharan considered the following: $A={\BC}[X,Y]$, $F=X^3+Y^3-1$ and $I=(F,T-1)\subset A[T]$.
Although $A$ is smooth, they produce a set of two generators of $I/(I^2T)$ that cannot be lifted to a set of two generators of $I$. This shows that Nori's question does not extend to two dimensional base rings in general, even when the ideal contains a monic. This example is based on the fact that $SK_1(A/F)$ is nontrivial. 

Since $SK_1$ of a curve over $\k$ is trivial, we wonder whether Nori's question can be answered in this case in the affirmative. In what follows we proceed to do so. It involves several steps. We begin with the following proposition due to R. G. Swan.

\bp \cite[Corollary 9.10]{R}\label{SP}
Let $A$ be a ring and $I$ be an ideal. Let $\gamma \in Sp_{2t}(A/I)$, $t \ge 1$. If
the class of $\gamma$ is trivial in $K_1Sp(A/I)$ and if $2t \ge \sr(A)-1$,
then $\gamma$ has a lift $\alpha \in Sp_{2t}(A)$.
\ep

\bl\label{l3.2}
Let $A$ be an affine domain of dimension $2$ over $\k$ and $\mathfrak{m}_1,\cdots,\mathfrak{m}_r$ be  maximal ideals of $A$. Let $S= A\smallsetminus (\mathfrak{m}_1\cup \cdots \cup \mathfrak{m}_r)$ and $I\subset S^{-1}A[T]$ be an ideal of height $2$. Then the natural map $\SL_2(S^{-1}A[T])\lra \SL_2(S^{-1}A[T]/I)$ is surjective. 
\el

\proof
There is an ideal $J\subset A[T]$ such that $S^{-1}J=I$. Let $K=\mathfrak{m}_1\cap \cdots \cap \mathfrak{m}_r$. There are two possibilities: $J+K[T]=A[T]$ or $J+K[T]\subsetneqq A[T]$.
In the first case, we have $\dim(S^{-1}A[T]/I)=0$ and we are done. 

Now we consider the second case when $J+K[T]$ is a proper ideal. In this case, it is easy to see that $\dim(A[T]/J)=\dim (S^{-1}A[T]/I)=1$. Note that $S^{-1}A[T]/I$ is the direct limit of affine $\k$-algebras of dimension one. It now follows from \cite[Corollary 2.4]{MKDIMRN} that $\SL_2(S^{-1}A[T])\lra \SL_2(S^{-1}A[T]/I)$ is surjective. \qed

\begin{proposition}\label{NQD2L}
	Let $A$ be an affine domain of dimension $2$ over $\k$ and $\mathfrak{m}_1,\cdots,\mathfrak{m}_r$ be some smooth maximal ideals of $A$. Let $S= A\smallsetminus (\mathfrak{m}_1\cup \cdots \cup \mathfrak{m}_r)$ and $R=S^{-1}A$. Let $I\subset R[T]$ be an ideal such that: (1) $I+\mathfrak{J}[T]=R[T]$, where $\mathfrak{J}$ is the Jacobson radical of $R$; (2) $\mu(I/I^2T)=\hh(I)=2$. Let $I=(f_1,f_2)+I^2T$ be given. Then there exist $F_1,F_2\in I$ such that $I=(F_1,F_2)$ and $F_i-f_i\in I^2T$ for $i=1,2$.  
\end{proposition}

\proof
By \cite{NMKL} there exists $e\in (I^2T)$ such that $I=(f_1,f_2,e)$ where $e(1-e)\in (f_1,f_2)$. Then, $I_e=R[T]_e=(1,0)$ and $I_{1-e}=(f_1,f_2)_{1-e}$. The unimodular row $(f_1,f_2)_{e(1-e)}$ can be completed to a matrix in $\SL_2(R[T]_{e(1-e)})$ and by a standard patching argument we obtain a surjection $P\sur I$ where $P$ is a projective $R[T]$-module of rank two with trivial determinant. As $R$ is smooth and semilocal, it follows that  $P$ is free and therefore, $I=(g_1,g_2)$. There exists a matrix $\overline{\sigma}\in GL_2(R[T]/I)$ such that $(\overline{f}_1,\overline{f}_2)=(\overline{g}_1,\overline{g}_2)\ol{\sigma}$. Let $\det(\ol{\sigma})=\ol{u}$ and let $\ol{u}\ol{v}=\ol{1}$ in $R[T]/I$. The unimodular row $(v,g_2,-g_1)\in \Um_3(R[T])$ is completable (as $R$ is smooth semilocal). Therefore, using \cite[Lemma 5.2]{SMBB3} we can find $h_1,h_2\in I$ such that $I=(h_1,h_2)$ and $(\overline{f}_1,\overline{f}_2)=(\overline{h}_1,\overline{h}_2)\ol{\theta}$ for some $\ol{\theta}\in \SL_2(R[T]/I)$. Applying Lemma \ref{l3.2} we can find a lift $\theta\in \SL_2(R[T])$ of $\ol{\theta}$. Let $(h_1,h_2)\theta=(H_1,H_2)$.

From the above paragraph, we have: $I=(H_1,H_2)$, where $H_i-f_i\in I^2$.
We still have to lift the generators of $I(0)$, namely, $f_1(0),f_2(0)$.
Since  $I+\mathfrak{J}R[T]=R[T]$, we have $I(0)=R$ and the rows $(f_1(0),f_2(0))$, 
$(H_1(0),H_2(0))$ are both unimodular. As $R$ is semilocal there is a matrix $\alpha\in E_2(R)$ such that $(f_1(0),f_2(0))= 
(H_1(0),H_2(0))\alpha$. Let $\alpha=\prod E_{ij}(a_{ij})$, $a_{ij}\in R$. As $I(0)=R$, there exists $\lambda_{ij}\in I$ such that $\lambda_{ij}(0)=a_{ij}$. Let $\Delta=\prod E_{ij}(\lambda_{ij})\in E_2(R[T])$. Taking $(F_1,F_2)=(H_1,H_2)\Delta$ we are done.
\qed

We shall recall some results which will be used to prove the main theorem in this section.

\bl\label{case2}
Let $A$ be a Noetherian domain and   $I \subset A[T]$ be an ideal such that $ I=\langle f_1,f_2\rangle+I^2T$. Write $J=I\cap A$ and assume further that $I_{1+J}=\langle g_1,g_2\rangle$ such that $g_i-f_i\in (I^2T)_{1+J}$ for $i=1,2$. Then there exist $F_1,F_2\in I$ such that 
$I=\langle F_1,F_2\rangle$ where 
$F_i-f_i\in I^2T$.
\el

\proof
We can clear denominators and assume that there is some $s\in J$ such that (abusing the notations) $I_{1+s}=\langle g_1,g_2\rangle$ and that
$g_i-f_i\in (I^2T)_{1+s}$ for $i=1,2$. Let $\alpha:A_{1+s}^2[T]\surj I_{1+s}$ be the corresponding surjection.

We have $I(0)=\langle f_1(0),f_2(0)\rangle$. Therefore, $A_s=I(0)_s=\langle f_1(0),f_2(0)\rangle$, implying further that
$I_s=A_s[T]=I(0)_s[T]=\langle f_1(0),f_2(0)\rangle$. Let $\beta: A_{1+s}^2[T]\surj I_s$ denote the corresponding surjection.

Let $K_1=\ker(\alpha_s)$ and $K_2=\ker(\beta_{1+s})$. Note that $(g_1,g_2)$ is a unimodular row over $A_{s(1+s)}[T]$. Since any unimodular row of length two is completable, $K_1$ is free of rank one. It is  easy to see that $K_2$ is also free of rank one.
Further, $\alpha_s(0)=\beta_{1+s}(0)$. By \cite[Lemma 2.9]{SMBB3}, there exists $\sigma\in GL_2(A_{s(1+s)}[T])$ such that $\sigma(0)$ is the identity matrix and $\beta_{1+s}\sigma=\alpha_s$. One can now use the well-known Quillen Splitting Lemma (see \cite[Lemma 2.10]{SMBB3}) and patch $\alpha$, $\beta$ to obtain a surjection $\phi:A^2[T]\surj I$ such that $\phi_{1+s}=\alpha$ and $\phi_s=\beta$. Now, $\phi$ yields $I=(F_1,F_2)$. It is now easy to check that these are the desired generators of $I$.
\qed

\bt \cite{AS} \label{SK1T}
Let $R$ be an affine algebra of dimension $1$ over $\k$. Then $SK_1(R)$ is trivial.
\et

Now we are ready to prove our main result in this section.
\bt\label{D2NQ}
Let $R$ be an affine domain of dimension two over $\k$. Let $I\subset R[T]$ be an ideal such that $\mu(I/I^2T)=\hh(I)=2$ and $R/(I\cap R)$ is smooth. Let $I=(f_1,f_2)+I^2T$ be given. Then there exist $F_1,F_2\in I$ such that $I=(F_1,F_2)$ and $F_i-f_i\in I^2T$ for $i=1,2$.
\et

\proof
Let $J=I\cap R$. Let `tilde' denote reduction modulo $(J^2T)$. We have $\widetilde{I}=(\widetilde{f}_1,\widetilde{f}_2)+\widetilde{(I^2T)}$. As $\dim(R[T]/(J^2T))\leq 2$, by Theorem \ref{IMMTT} there exist $g_1,g_2\in I$ such that $\widetilde{I}=(\widetilde{g}_1,\widetilde{g}_2)$ such that $\widetilde{g}_i-\widetilde{f}_i\in\widetilde{(I^2T)}$. Therefore,
$I=(g_1,g_2)+(J^2T)$ such that $g_i-f_i\in (I^2T)$. Using \cite{NMKL} there exists $e\in J^2T$ such that $I=(g_1,g_2,e)$ and $e(1-e)\in (g_1,g_2)$. Moreover by a theorem due to Eisenbud-Evans (\cite{EE} also see, \cite[Corollary 2.13]{SMBB3} ) replacing $g_i$ by $g_i+\lambda_ie$ (and retaining the same notations) we may assume that $\hh((g_1,g_2)_e)\ge 2$. Let $I'=(g_1,g_2,1-e)$. Then we have  
$I'+(J^2T)=R[T]$, $\hh(I')\geq 2$ and $I\cap I'=(g_1,g_2)$.

If $I'=R[T]$, then we are done. Therefore, we assume that $I'$ is proper and  $\hh(I')=2$. 
We have $I'=(g_1,g_2)+I'^2$. Note that $I'(0)=R$. Applying \cite[Remark 3.9]{BR} we can lift $g_1,g_2$ so that $I'=(h_1,h_2)+(I'^2T)$ where $h_i-g_i\in I'^2$ for $i=1,2$.

Let $J'=I'\cap R$. Let $B=R_{1+J}$ and $C=B_{1+J'}=R_{1+J+J'}$.
Note that since $R/J$ is smooth, the ideal of singular locus of $R/J$ (which is extended from $R$) is co-maximal with $J$, hence the ring $B$ is smooth. This implies that $C$ is smooth, being further localization of a smooth ring. It has been proved in \cite[Theorem 3.8]{BR} that the ring $C$ is semilocal.

We have $I'C[T]=(h_1,h_2)+(I'^2T)$. Since $I'C[T]+(J^2T)C[T]=C[T]$ and $J$ is contained in the Jacobson radical of $C$, we can apply Theorem \ref{NQD2L} and ensure that $h_1,h_2$ can be lifted to a set of generators of $I'C[T]$. Now,
we can apply Lemma \ref{PLBR} and obtain:
$$I'B[T]=(k_1,k_2)\text{ such that } k_i-h_i\in (I'^2T)B[T].$$

Note that, in view of Lemma \ref{case2}, to prove the theorem it will be enough to show that $IB[T]=(\ga_1,\ga_2)$ such that $\ga_i-g_i\in (I^2T)B[T]$. The remaining part of the proof is dedicated to show this.

We have $I'B[T]+(J^2T)B[T]=B[T]$. Let us write $D=B[T]/J^2B[T]$ and
`bar' denote modulo $J^2B[T]$.  Now, $(\overline{k_1},\ol{k_2})\in \Um_2(D)$. As $ (\overline{k_1},\ol{k_2})$ is a unimodular row of length two, there is a matrix $\sigma\in \SL_2(D)$ such that $(\overline{k_1},\ol{k_2})\sigma=(\ol{1},\ol{0})$. 

\smallskip

\noindent
{\bf Claim:} $\sigma$ can be lifted to a matrix $\tau\in SL_2(B[T])$.

\noindent
\emph{Proof of the claim.} Since $B[T]=R_{1+J}[T]$, we observe that $B$ is the direct limit of affine $\k$-algebras of dimension two and therefore, applying \cite[Corollary 17.3]{SV} we obtain: $\sr(B[T])\leq max\{2, \dd(B[T])\}=3$.

Let us now consider $SK_1(D)$ and $K_1Sp(D)$. We have $$D_{\text{red}}=B[T]/\sqrt{J}[T]=B[T]/J[T]=(R/J)[T],$$
since $J$ is reduced. Since $R/J$ is smooth, we have $SK_1((R/J)[T])=SK_1(R/J)$ and $K_1Sp((R/J)[T])=K_1Sp(R/J)$.
Since $\dd(R/J)=1$, by Theorem \ref{SK1T} we have $SK_1(R/J)=0$ and applying \cite[Lemma 16.2 and the remark following the lemma]{SV} we further obtain that $K_1Sp(R/J)=SK_1(R/J)$ and hence it is trivial as well.

For any ring $S$, by a result of Bass \cite[IX, 1.3]{bass}, we have $SK_1(S)=SK_1(S_{\text{red}})$. We conclude that $SK_1(D)=0$. On the other hand, for a ring $C$, by \cite[9.9, 9.12]{R} the natural map $K_1Sp(C)\lra K_1Sp(C_{\text{red}})$ is injective. Therefore, from the above computation we see that
$K_1Sp(D)=0$.

We can now apply Swan's result (Proposition \ref{SP}) with $t=1$. Since $Sp_2(B[T])$ is the same as $\SL_2(B[T])$, the claim is proved.

\smallskip

Let $(k_1,k_2)\tau=(k'_1,k'_2)$. Then $I'=(k'_1,k'_2)$ and $k'_1\equiv 1$ modulo $J^2B[T]$ and $k'_2\equiv 0$ modulo $J^2B[T]$. Write $(\gb_1,\gb_2)=(k'_1,k'_2)\epsilon$, where $\epsilon=\begin{pmatrix}
1 & 1\\
0 & 1
\end{pmatrix}$, (here $\gb_1=k'_1$ and $\gb_2=k'_1+k'_2$). Then $\gb_i\equiv 1$ modulo $J^2B[T]$ for $i=1,2$. 

We now write $B[T]=A$ and introduce a new variable $X$ and consider the following ideals in $A[X]$:

$$K'=(\gb_1,X+\gb_2),\,\,\, K''=IA[X],\,\,\,  K=K'\cap K''$$ 
Let us write the ideal $\gb_1 A$ as $\mathfrak{n}$.
We have $K_{1+\mathfrak{n}}=K'_{1+\mathfrak{n}}=(\gb_1,X+\gb_2)$. Recall that we have $I\cap I'=(g_1,g_2)$, implying that $(I\cap I')A=(g_1,g_2)A$. Let $(g_1,g_2)\tau=(g'_1,g'_2)$ and write
$(l_1,l_2)=(g'_1,g'_2)\epsilon$. Then also we have 
$(I\cap I')A=(l_1,l_2)A$ and all the relations are retained.

Now $  (l_1,l_2)A_{1+\mathfrak{n}}=I'A_{1+\mathfrak{n}}=
K'_{1+\mathfrak{n}}(X=0)=K_{1+\mathfrak{n}}(X=0)$.

We also have $K_{1+\mathfrak{n}}(X=0)=(\gb_1,\gb_2)$.
Then note that $(l_1-\gb_1,l_2-\gb_2)=(g_1-k_1,g_2-k_2)\tau \epsilon\in I'^2A_{1+\mathfrak{n}}\times I'^2A_{1+\mathfrak{n}}$ as $g_i-k_i=(g_i-h_i)+(h_i-k_i)\in (I'^2B[T])$.
Since $K'_{1+\mathfrak{n}}(X=0)=(\gb_1,\gb_2)=(l_1,l_2)$ and two sets of generators of $K'_{1+\mathfrak{n}}(X=0)$ differ by an element of $\GL_2(A_{1+\mathfrak{n}} )$, there exists $\alpha\in \GL_2(A_{1+\mathfrak{n}} )$ such that $(\gb_1,\gb_2)\alpha=(l_1,l_2)$. Let $(\gb_1,X+\gb_2)\alpha=(G_1(X),G_2(X))$, then $G_i(0)=l_i$ for $i=1,2$.

Recall that $\mathfrak{n}=(\gb_1)$ is comaximal with $J^2B[T]$. We choose some $s\in \mathfrak{n}$ such that $1+s\in J^2B[T]$ and 
$K'A_{1+sA}=(G_1(X),G_2(X))$ with $G_i(0)=l_i$ for $i=1,2$. 

Let $\phi:A_{1+sA}[X]^2\sur K_{1+sA}$ be the surjection corresponding to $K_{1+sA}=(G_1(X),G_2(X))$. We have a surjection $\psi: A_{s}[X]^2\sur K_s$ induced by the following: $K_s=K''_s=I_s=(l_1,l_2)$.

The surjections $\phi_s: A_{s(1+sA)}[X]^2\sur K_{s(1+sA)}$ and $\psi_{1+sA}:A_{s(1+sA)}[X]^2\sur K_{s(1+sA)}$ agree when $X=0$. As both the kernels are free applying \cite[Lemma 2.17]{MPHIL} there exists $\kappa(X)\in \SL_2(A_{s(1+sA)}[X])$ such that $\phi_s\kappa(X)=\psi_{1+sA}$ and $\kappa(0)=\text{Id}$. Hence by Quillen's splitting lemma \cite{Q} $\kappa(X)$ splits. Therefore, from the universal property of the fiber product we obtain that $K=(H_1(X),H_2(X))$ such that $H_i(0)=l_i$ for $i=1,2$.

Now, $I=K(1-\gb_2)=(H_1(1-\gb_2), H_2(1-\gb_2))$. We write $H_i(1-\gb_2)=\ga_i$.  As the constant term of $H_i$ is $l_i$ and $\gb_2\equiv 1$ modulo  $J^2B[T]$, it follows that $IA=(\ga_1,\ga_2)$ with $\ga_i- l_i$ modulo $J^2B[T]$.

Let us now revert back to the original notations (recall: $A=B[T]=R_{1+J}[T]$).
We have thus far been able to establish the following:
$$I_{1+J}=(\ga_1,\ga_2)\text{ such that } \ga_i- l_i\in  J^2B[T]$$
Also recall that we started with  $I=(g_1,g_2)+(J^2T)$ and then applied some automorphisms on $(g_1,g_2)$ to get $(l_1,l_2)$. So we have $I=(l_1,l_2)+(J^2T)$. However $\ga_1,\ga_2$ are lifts modulo $J^2[T]$. We need to find $\gamma_1,\gamma_2$ so that
$I_{1+J}=(\gamma_1,\gamma_2)$ with $\gamma_i-l_i\in J^2[T]_{1+J}$ and $\gamma_i(0)=l_i(0)$ for $i=1,2$. Once we have done this we can apply inverses of the said automorphisms on $(\gamma_1,\gamma_2)$ to solve the problem. The remaining part of the proof is dedicated to find such $\gamma_i$.

We have $I(0)_{1+J}=(l_1(0),l_2(0))=(\ga_1(0),\ga_2(0))$ such that $\ga_i(0)-l_i(0)\in (J^2)_{1+J}$. 

Note that $J\subset I(0)$ and $J^2\subset JI(0)$. Therefore, we can write
$\ga_1(0)-l_1(0)=c\ga_1(0)+d\ga_2(0)$, where $c,d\in JB$. Similarly,
$\ga_2(0)-l_2(0)=e\ga_1(0)+f\ga_2(0)$, where $e,f\in JB$. Putting it in another way,
$$(\ga_1(0),\ga_2(0))\delta
=(l_1(0),l_2(0))$$
where $\delta=\begin{pmatrix}

1-c & -e\\

-d & 1-f
\end{pmatrix}$.
Note that the determinant of the above matrix is $1$ modulo $J$. Since $J$ is contained in the Jacobson radical of $R_{1+J}$, it is an invertible matrix in $R_{1+J}$.

Let $(\gamma_1,\gamma_2)= (\ga_1,\ga_2)\delta$. 
Then note that:
\begin{enumerate}
\item $I_{1+J}=(\gamma_1,\gamma_2)$, as $\delta\in \GL_2(R_{1+J})$;
\item $\gamma_i-l_i\in I^2_{1+J}$, as $\gamma_1-l_1=(\gamma_1-\alpha_1)-(\alpha_1-l_1)=-c\alpha_1-d\alpha_2+(\alpha_1-l_1)\in I^2_{1+J}$ and $\gamma_2-l_2=(\gamma_2-\alpha_2)-(\alpha_2-l_2)=-e\alpha_1-f\alpha_2+(\alpha_2-l_2)\in I^2_{1+J}$; 
\item  $\gamma_i(0)=l_i(0)$ for $i=1,2$.
\end{enumerate} 
This completes the proof.\qed

\section{Addressing the failures: Precise obstruction}
As mentioned earlier, in the appendix of a paper by Mandal \cite{SM}, Nori asked the following question, which is motivated by certain results in topology. For the convenience of reading, we restate the ``free" version of the question below.
\smallskip

\quest
Let $R$ be a smooth affine domain of dimension $d$ over an infinite perfect field $k$ and $I\subset R[T]$ be an ideal of height $n$ such that $\mu(I/I^2T)=n$, where $2n\geq d+3$. Assume that $I=(f_1,\cdots,f_n)+(I^2T)$ is given. Then, do there exist $F_i\in I$ ($i=1,\cdots,n$) such that 
$I=(F_1,\cdots,F_n)$ where $F_i-f_i\in (I^2T)$ for $i=1,\cdots,n$?
\smallskip

If $I$ contains a monic polynomial, then S. Mandal \cite{SM} proved that the answer is in the affirmative where he needs the ring $R$ to be just Noetherian. Ideals containing monic polynomials are of a different league and let us leave them out of our discussion. So, from now on we assume that $I$ does not contain monic in the above question. Nori's question has been answered comprehensively. First, Mandal-Varma \cite{mv} proved it to be true when $R$ is local. Bhatwadekar-Sridharan \cite{BR} gave an affirmative answer when $n=d\geq 3$ and $k$ is infinite perfect. Later, Bhatwadekar-Keshari settled it in the affirmative \cite{BHMK}  for $2n\geq d+3$ with the same assumption on $k$. 

On the other hand, Bhatwadekar-Mohan Kumar-Srinivas gave an example in \cite[Example 6.4]{BR} to show that Nori's question will have a negative answer if $R$ is not smooth (even when $R$ is local). They constructed an example of a normal affine $\mathbb C$-domain $R$ of dimension 3 which has an isolated singularity at the origin and an ideal $I\subset R[T]$ of height 3 such that a given set of generators of $I/(I^2T)$ cannot be lifted to a set of generators of $I$. 

The results and the example stated above had profound impact on the development of the theory in understanding the behaviour of projective modules and local complete intersection ideals in past twenty years. Among recent instances, the Bhatwadekar-Sridharan solution played a crucial role in computing the group of isomorphism classes of oriented stably free $R$-modules of rank $d$ where $R$
is a smooth affine domain of dimension $d$ over $\br$ (\cite{DTZ}, see also \cite{DTZ2}). Further, Asok-Fasel \cite{af} used it successfully to establish the isomorphism between the $d$-th Euler class group  and the $d$-th Chow-Witt group (also the isomorphism between the weak Euler class group and the Chow group) --- thus establishing a long standing conjecture.  
\smallskip

In this context, we delve deep into this phenomenon and observe that in the example of Bhatwadekar-Mohan Kumar-Srinivas \cite[Example 6.4]{BR}, the ideal $I\subset R[T]$ is such that: $R[T]/I$ is smooth but $I\cap R$ is contained in a non-smooth maximal ideal (see the local version of their counter-example for more clarity). We now pose the following rephrased question.

\quest
Let $R$ be  an affine domain of dimension $d$ over an infinite perfect field $k$ and $I\subset R[T]$ be an ideal of height $n$ such that $\mu(I/I^2T)=n$, where $2n\geq d+3$. Assume that $I=(f_1,\cdots,f_n)+(I^2T)$ is given. Then, what is the precise obstruction for $I$ to have a set of generators   $F_1,\cdots,F_n$ such that 
$F_i-f_i\in (I^2T)$ for $i=1,\cdots,n$?
\medskip

Obviously we have left out the case when $I$ contains a monic polynomial. We prove that the obstruction lies in the fact as to whether
$I\cap R$ is contained in only smooth maximal ideals or not. More precisely, we prove the following result. We have decided to give the details for the case $n=d$. We shall comment on the other versions in the sequel.

\bt\label{free}
Let $R$ be an affine domain of dimension $d$ over an infinite perfect  field $k$ and $I\subset R[T]$ be an ideal of height $n$ such that  $\mu(I/I^2T)=n$, where $2n\geq d+3$. Assume that $I=(f_1,\cdots,f_n)+I^2T$ is given. Assume further that $I\cap R$ is contained only in smooth maximal ideals of $R$. Then, there exist $F_i\in I$ such that 
$I=(F_1,\cdots,F_n)$, where $F_i-f_i\in I^2T$, $i=1,\cdots,n$. 
\et

To prove the above result, the most crucial proposition is the following improvement of \cite[Theorem 3.8]{BR}. Note that we do not assume $P$ to be extended.

\bp\label{patch}
Let $A$ be a  domain which is essentially of finite type over a field. Let $I\subset A[T]$ be an ideal such that $J:=I\cap A$ is contained only in smooth maximal ideals. Let $P$ be a projective $A[T]$-module such that there is a surjection
$$\ol{\varphi}:P\sur I/(I^2T).$$ Assume that there is a surjection
$$\theta: P_{1+J}\sur I_{1+J}$$such that $\theta$ is a lift of $\ol{\varphi}\otimes A_{1+J}$. Then, there is a surjection $\Phi:P\sur I$ which lifts $\ol{\varphi}$.
\ep

\proof
From the map $\theta$, clearing denominators we can find $s_1\in J$ such that $\theta:P_{1+s_1}\sur I_{1+s_1}$ is surjective (we are using the same notation $\theta$). 

Let $A=S^{-1}B$, where $B$ is an affine $k$-domain. Since singular locus of $B$ is a closed set $V(\MI)$, singular locus of $A$ is 
$V(S^{-1}\MI)$. Now $A_{1+J}$ is regular, so $J+S^{-1}\MI=A$ and so $1+s_2\in S^{-1}\MI$ for some $s_2\in J$ and $A_{1+s_2}$ is regular.

We now take $(1+s):=(1+s_1)(1+s_2)$ and consider $\theta:P_{1+s}\sur I_{1+s}$. Note that, as $A_{1+s}$ is a regular ring containing a field, by a result of Popescu \cite{pop}, the module $P_{1+s}$ is extended from $A_{1+s}$. 

The map $\ol{\varphi}$ induces a surjection, say, $\ol{\varphi}(0):P/TP\sur I(0)$. As $s\in J$, we have $I(0)_s[T]=I_s=A_s[T]$. Therefore,
we have $\ol{\varphi}(0)_s :(P/TP)_s\sur I_s$. Then we have a surjection $\ga:P_s\sur I_s$ (composing  $\ol{\varphi}(0)_s $ with the canonical map $P_s\sur (P/TP)_s$).

We now proceed to patch the two maps $\theta:P_{1+s}\sur I_{1+s}$ and $\ga:P_s\sur I_s$. We move to $A_{s(1+s)}[T]$. As $P_{s(1+s)}$ is extended from $A_{s(1+s)}$, we have a projective $A_{s(1+s)}$-module, say, $P'$ such that $P'[T]=P_{s(1+s)}$. We finally have:
$$\theta_s:P'[T]\sur I_{s(1+s)} \,(=A_{s(1+s)}[T]), \text{ and }$$
$$\ga_{1+s}:P'[T]\sur I_{s(1+s)} \,(=A_{s(1+s)}[T]),$$
where $\theta_s$ and $\ga_{1+s}$ are equal modulo $(T)$. Since $A_{s(1+s)}$ is a regular ring containing a field, we also note that the kernels of $\theta_s$ and $\ga_{1+s}$ are both extended from $A_{s(1+s)}$. Therefore, by \cite[Lemma 2.9]{SMBB3} there is an isomorphism $\sigma:P'[T]\iso P'[T]$ such that $\sigma(0)=id$ and $\ga_{1+s}\sigma=\theta_s$. We can now patch $\theta:P_{1+s}\sur I_{1+s}$ and $\ga:P_s\sur I_s$ using Plumstead's patching technique (see \cite{P}) to obtain a surjection $\Phi:P\sur I$. It is then easy to check that $\Phi$  lifts $\ol{\varphi}$.
\qed

We now present the following \emph{``projective"} version of Theorem \ref{free} mentioned above. This is an improvement  of 
the result of Bhatwadekar-Sridharan \cite[Theorem 3.8]{BR}. The proof is essentially contained in \cite{BR}. We just give a sketch and for the details we refer to their paper.

\bt\label{proj}
Let $R$ be  a domain of dimension $d$ which is essentially of finite type over an infinite perfect field $k$. Let $n$ be an integer such that $2n\geq d+3$. Let $I\subset R[T]$ be an ideal of height $n$ such that $J:=I\cap R$ is contained only in smooth maximal ideals. Let $P$ be a projective $R[T]$-module of rank $n$ such that there is a surjection
$$\ol{\varphi}:P\sur I/(I^2T).$$ Then, there is a surjection $\Phi:P\sur I$ which lifts $\ol{\varphi}$.
\et

\proof
Let us write $A=R_{1+J}$. Then $A$ is smooth and essentially of finite type over $k$. We can therefore apply \cite[Theorem 4.13]{BHMK} on $A$, $I_{1+J}$, and $P_{1+J}$ to obtain that there is a surjection $\psi:P_{1+J}\sur I$ that is a lift of $\varphi\otimes R_{1+J}$.
Now  Proposition \ref{patch} implies that we have a surjetion $\Phi:P\sur I$ which lifts $\ol{\varphi}$.
\qed

\rmk
The above theorem shows that the condition that $I\cap R$ is contained only smooth maximal ideals is sufficient to find surjective 
lifts, as mentioned in the questions above.  

\rmk
In the above theorems we can remove the condition that $k$ is perfect. When $R$ is local, this was done in \cite{Das2012}. To extend it to the semilocal version, one has to use an induction argument on the number of maximal ideals, as done in \cite{MKD1}. We leave these details to the readers to check. 

\rmk
One may wonder whether instead of the smoothness condition on $I\cap R$, we can impose it on $I$ itself. Unfortunately, that would not work. In the example of Bhatwadekar-Mohan Kumar-Srinivas mentioned above, $R[T]/I$ is smooth and the lifting fails.

\bt\label{monic}
Let $R$ be a domain of dimension $d$ containing a field $k$ (no restriction on $k$) and $n$ be an integer such that $2n\geq d+3$.  Let $I\subset R[T]$ be an ideal of height $n$ such that $J:=I\cap R$ is contained only in smooth maximal ideals. Let $P$ be a projective $R[T]$-module of rank $n$ and $\varphi:P\sur I/I^2T$ be a surjection. Assume that there is a surjection $\varphi':P\otimes R(T)\sur IR(T)$ such that $\varphi'$ is a lift of $\varphi\otimes R(T)$. Then there is a surjection $\psi:P\sur I$ which lifts $\varphi$.
\et

\proof Let us take $J=I\cap R$. First we observe that, in view of Proposition \ref{patch} it is enough to prove the theorem on $R_{1+J}$. Moreover, one may note that $R_{1+J}$ is a regular domain containing $k$. Therefore, the proof follows from \cite[Proposition 4.9]{BHMK}.\qed

We now proceed to show some interesting applications of Proposition \ref{patch} and Theorem \ref{proj}. For a $\bq$-algebra $R$ of dimension $d\geq 3$, the $d$-th Euler class group $E^d(R[T])$ was defined in \cite{MKD1}. It was further proved that the canonical map $\phi:E^d(R)\lra E^d(R[T])$ is injective. The morphism $\phi$ is an isomorphism if $R$ is smooth  but it may not be surjective if $R$ is not smooth (see \cite{MKD1} for the details). In this context, we may ask, precisely which Euler cycles $(I,\omega_I)\in E^d(R[T])$ have a preimage in $E^d(R)$? We answer this question in the following form.

\bt\label{ECIT}
Let $R$ be an affine domain of dimension $d\geq 3$ over a field $k$ of characteristic zero. Let $(I,\omega_I)\in E^d(R[T])$ be such that
$I\cap R$ is contained only in smooth maximal ideals. Then  $(I,\omega_I)$
is in the image of the canonical morphism $\phi: E^d(R)\lra E^d(R[T])$.  
\et

\proof

As $k$ is infinite, applying \cite[Lemma 3.3]{BR}  we can find $\sigma\in k$ such that either $I(\sigma)=R$ or $\hh(I(\sigma))=d$. Changing $T$ by $T-\sigma$, we may assume that either $I(0)=R$ or $\hh(I(0))=d$. If $I(0)=R$, then $\omega_I$ can be lifted to a surjection $\alpha:R[T]^d\sur I/(I^2T)$. Then $\alpha$ lifts to a surjection from $R[T]^d$ to $I$ and consequently, $(I,\omega_I)=0$. Therefore, we assume that $I(0)$ is proper of height $d$. Then $(I,\omega_I)$ induces $(I(0),\omega_{I(0)})\in E^d(R)$. If $(I(0),\omega_{I(0)})=0$ in $E^d(R)$, then also $\omega_I$ can be lifted to a surjection $\alpha:R[T]^d\sur I/(I^2T)$ and we will be done by Theorem \ref{proj} (taking $P$ to be free). So let $(I(0),\omega_{I(0)})\neq 0$ in $E^d(R)$.

Using the moving lemma \cite[Corollary 2.14]{SMBB3} together with Swan's Bertini Theorem \cite[Theorem 2.11]{BRS}, we can find a reduced  ideal $K\subset R$ of height $d$ which is comaximal with $I\cap R$ and a local orientation $\omega_K$ such that $(I(0),  \omega_{I(0)})+(K,\omega_K)=0$ in $E^d(R)$.

As $K$ is reduced of height $d$ and $\mu(K/K^2)=d$, we observe that $K$ is product of a finite number of distinct smooth maximal ideals of $R$.

Let $L=I\cap K[T]$. The local orientations $\omega_I$ and $\omega_K$ will induce $\omega_L$ and we have 
$$(L,\omega_L)= (I,\omega_I)+(K[T],\omega_{K[T]}) \text{ in } E^d(R[T]).$$

As $(L(0),\omega_{L(0)})= (I(0),\omega_{I(0)})+(K,\omega_{K})=0$, it follows that $\omega_L$ can be lifted to a surjection 
$\lambda:R[T]^d\sur L/(L^2T)$.

Now $L\cap R=(I\cap R)\cap K$. Since $K$ is reduced and $I\cap R$ is contained only in smooth maximal ideals of $R$, it follows that $L\cap R$ is contained only in smooth maximal ideals of $R$. Therefore, by Theorem \ref{proj} $\lambda$ can be lifted to a surjection $\alpha:R[T]^d\sur L$. As a consequence, $(L,\omega_L)=0$ in $E^d(R[T])$ and we have 
$$(I,\omega_I)=-(K[T],\omega_{K[T]})\in \phi (E^d(R)).$$
\qed

\rmk
With notations as above, let $(I,\omega_I)\in E^d(R[T])$ be such that $I$ is a non-extended ideal of $R[T]$ and $I$ does not contain a monic polynomial.
If we further drop the condition that  $I\cap R$ is contained only in smooth maximal ideals, then there is an example \cite[Remark 3.4]{d2} which shows that $(I,\omega_I)$ may not be in the image of $\phi$.

\medskip

\rmk
The couple of theorems given above were proved in \cite[Proposition 4.9 and Theorem 4.13]{BHMK}  assuming $R$ to be smooth. We remark that \cite[Proposition 4.9]{BHMK} is crucially used to prove \cite[Theorem 4.13]{BHMK}.

Let $A$ be a commutative Noetherian ring and $P$ be a projective $A[T]$-module. A necessary condition for $P$ to have a unimodular element is that the $A$-module $P/TP$ should have a unimodular element. Based on our ideas above, we prove the following result by imposing a sufficient condition on a \emph{generic section} of $P$.

\bt\label{unimod}
Let $R$ be  an affine domain of dimension $d\geq 3$ over an infinite  field $k$ and $P$ be a projective $R[T]$-module of rank $d$ with trivial determinant. Assume that $P/TP$ has a unimodular element and that there is a surjection $\alpha:P\sur I$ where $I\subset R[T]$ is an  ideal of height $d$ such that $J:=I\cap R$ is contained only in smooth maximal ideals. Then $P$ has a unimodular element.
\et 

\proof
As $k$ is infinite, there is some $\lambda\in k$ such that $I(\lambda)=R$ or $\hh(I(\lambda))=d$. We note that the automorphism $T\mapsto T-\lambda$ preserves all the hypotheses of the theorem. Therefore, we may assume that $\hh(I(0))=d$.

Let us fix an isomorphism $\chi:R[T]\iso \wedge^d (P)$. We note that $P/IP$ is free and we choose an isomorphism 
$\sigma:(R[T]/I)^d\iso P/IP$ such that $\wedge^d \sigma=\chi\otimes R[T]/I$. Writing $\overline{\alpha}$ for $\alpha\otimes R[T]/I$ we then have:
$$(R[T]/I)^d \stackrel{\sigma}{\iso} P/IP\stackrel{\overline{\alpha}}{\sur} I/I^2$$
Let us call this composite $\overline{\alpha}\sigma$ as $\omega_I:(R[T]/I)^d\sur I/I^2$.

Let $\alpha(0)$ be the surjection from $P/TP$ to $I(0)$ be the one induced by $\alpha$. Putting $T=0$ in everything above we obtain 
$\omega_{I(0)}:(R/I(0))^d\sur I(0)/I(0)^2$, which is easily seen to be the one induced by $\omega_I$ by assigning $T=0$.

As $P/TP$ has a unimodular element, it follows from \cite{SMBB3} that $\omega_{I(0)}$ can be lifted to a surjection $\theta:R^d\sur I(0)$. Now, applying \cite[Remark 3.9]{BR} we see that $\omega_I$ can be lifted to a surjection $\beta:(R[T]/I)^d\sur I/(I^2T)$. As $I\cap R$ is contained only in smooth maximal ideals, it follows from Theorem \ref{proj} that $\beta$ can be lifted to a surjection $\gamma: R[T]^d\sur I$. Note that $\gamma$ also lifts  $\omega_I$.

We now move to the ring $R(T)$, which is obtained from $R[T]$ by inverting all the monic polynomials. We see that 
$(P,\chi\otimes R(T),\alpha\otimes R(T))$ induces $\omega_I\otimes R(T)$ and $\omega_I\otimes R(T)$ has a surjective lift, namely, $\gamma\otimes R(T)$. 
It then follows from \cite[Corollary 3.4]{SMBB3} that $P\otimes R(T)$ has a unimodular element. Since $P/TP$ has a unimodular element, we are done by \cite[Theorem 5.2 and Remark 5.3]{SMBHLRR}.
\qed

\smallskip

\rmk
The above result can also be proved when the determinant of $P$ is not trivial. We briefly explain the arguments. One can easily notice that the above proof also works when the determinant of $P$ is extended from $R$. In the case when the determinant is not extended, we can use the arguments given by Bhatwadekar in \cite[Proposition 3.3]{SMBmonic} and use \cite[Lemma 3.2]{SMBmonic} to 
reduce to the case when the determinant is extended.

\medskip

Now we proceed to prove  a result on projective generation of curves. Let us first digress a little. Let $A$ be a ring of dimension $d$ and $L\subset A$ be an ideal of height $d$ such that $\mu(L/L^2)=d$. Let $P$ be a projective $A$-module of rank $d$ with trivial determinant. Let $\chi:A\iso \wedge^d P$ be an isomorphism. Assume that there is a surjection $\theta:P\sur L$. As $P/LP$ is a free $A/L$-module of rank $d$, we can choose an isomorphism $\sigma:(A/L)^d\iso P/LP$ such that $\wedge^d \sigma=\chi\ot A/L$ and obtain (writing $\ol{\theta}:=\theta\ot A/L$):
$$(A/L)^d\stackrel{\sigma}{\iso} P/LP\stackrel{\overline{\theta}}{\sur} L/L^2.$$
Let us call $\omega_L:=\overline{\theta}\sigma:(A/L)^d\sur L/L^2$. In this situation, we say that $\omega_L$ is obtained from the triple $(P,\chi,\theta)$.

\bt\label{projgen}
Let $R$ be  an affine domain of dimension $d\geq 2$ over an infinite  field $k$ and let $I\subset R[T]$ be an ideal of height $d$ such that: (i) $I/I^2$ is generated by $d$ elements; (ii) $J:=I\cap R$ is contained only in smooth maximal ideals. Let $\lambda\in k$ be such that $\hh(I(\lambda))\geq d$. Let $\omega_I:(R[T]/I)^d\sur I/I^2$ be a surjection such that the induced surjection $\omega_{I(\lambda)}:(R/I(\lambda))^d\sur I(\lambda)/I(\lambda)^2$ is obtained from a triple $(P,\chi,\theta)$ (as described above).
Then  $I$ is surjective image of  $P[T]$. 
\et

\proof
First of all, for any commutative Noetherian ring $B$ and for any ideal $L\subset B$ such that $\mu(L/L^2)=2$, there is a projective
 $B$-module $P$  mapping onto $L$. One does not need any of the hypotheses as above. Therefore, we assume that $d\geq 3$ from now on.

Next, changing $T$ by $T-\lambda$, we may assume $\lambda=0$.

Further, if $I(0)=R$, then by \cite[Remark 3.9]{BR} we can lift $\omega_I$ to a surjection $\omega_{I}':R[T]^d\sur I/(I^2T)$. Then we can apply Theorem \ref{proj} to conclude that $I$ is generated by $d$ elements. Therefore, we assume that $I(0)$ is a proper ideal of height $d$.

By hypothesis, there is a projective $R$-module $P$ of rank $d$, an isomorphism $\chi:R\iso \wedge^d P$ and a surjection $\theta:P\sur I(0)$ such that $\omega_{I(0)}$ is the composite (writing $\overline{\theta}:=\theta\otimes R/I(0)$):
$$(R/I(0))^d\stackrel{\sigma}{\iso} P/I(0)P\stackrel{\overline{\theta}}{\sur} I(0)/I(0)^2.$$

Now $\chi(T):R[T]\iso \wedge^d P[T]$ is an isomorphism. As $P[T]/IP[T]$ is a free $R[T]/I$-module, we can choose an isomorphism 
$\tau:(R[T]/I)^d\iso P[T]/IP[T]$ such that $\wedge^d\tau =\chi(T)\ot R[T]/I$. Note that $\tau(0)$ and $\sigma$ differ by a matrix, say, $\beta \in SL_d(R/I(0))$. As $\dd(R/I(0))=0$, we have $SL_d(R/I(0))=E_d(R/I(0))$ and therefore we can lift $\beta$ to $\E_d(R[T]/I)$ and adjust $\tau$ so that $\tau(0)=\sigma$.

We have the composite surjection
$$P[T] \stackrel{p}{\sur} P[T]/IP[T]\stackrel{\tau^{-1}}{\iso} (R[T]/I)^d\stackrel{\omega_I}{\sur} I/I^2,$$
where $p$ is the natural projection map. 
Let us call $\alpha:=\omega_I\circ \tau^{-1}\circ p$.

As $\tau(0)=\sigma$, it is clear that we have: $\alpha:P[T]\sur I/I^2$ and $\theta:P\sur I(0)$ such that $\alpha(0)=\theta\ot R/I(0) (=\ol{\theta})$. Now, applying \cite[Remark 3.9]{BR} we can lift $ \alpha$ to a surjection $\delta:P[T]\sur I/(I^2T)$. Applying Theorem \ref{proj} we further lift $\delta$ to a surjection $\varphi:P[T]\sur I$.
\qed

\rmk
The reader may compare the above result with \cite[Theorem 3.5]{MKDSMB}, where the condition on $I\cap R$ was not assumed but $k$ was assumed to be of characteristic zero.

\section{Improvement I: A theorem of Mandal}
In this section we focus on a question asked by M. V. Nori (Theorem \ref{1}) and its subsequent developments. We improve the bound imposed by S. Mandal \cite{SM}, for non extended ideals of affine algebras over $\ol{\mathbb{F}}_p$, which contains a monic polynomial. Before going to our main theorem we shall state the following lemmas. The proof of the following lemma is standard \cite[see Lemma 2 and Proposition 2]{P}, hence we opt to omit the proof

\bl\label{NQS}
Let $A$ be a commutative Noetherian ring and $I\subset A[T]$ be any ideal containing a monic polynomial and $J=I\cap A$. Suppose that $P$ is a projective $A$-module of rank $n\ge 2$ and $\ol{\phi}:P[T]\surj I/(I^2T)$ is a surjection. Moreover, assume that there exists $j\in J^2$ and a surjective map $\phi'(T):P_{1+j}[T]\surj I_{1+j}$, which lifts $\ol{\phi}\otimes A_{1+j}[T]$. Then there exists a surjective map $\phi(T):P[T]\surj I$, which lifts $\ol{\phi}$.
\el
%

\bt \label{1}
Let $A$ be an affine algebra over  $\overline{\mathbb{F}}_p$ and $I\subset A[T]$ be an ideal containing a monic polynomial. Assume that
$I=(f_1,...,f_n)+I^2T$, where  $n\ge \max\{(\dim\frac{A[T]}{I}+1),2\}$. Then there exist $F_i\in I$ such that
$I=(F_1,...,F_n)$ with $F_i-f_i\in I^2T$ for all $i$.
\et 
\proof We shall divide the proof into the following two cases:
\paragraph{\textbf{Case - 1}}
In this case we shall assume that $n\ge 3$. Let $f\in I$ be a monic polynomial. Without loss of generality we may assume that $f_1$ is monic (by replacing $f_1$ with $f_1+T^kf^2$, for some suitably chosen $k>0$). By \cite{NMKL} there exists $e\in (I^2T)$ such that $I=(f_1,...,f_n,e)$ with $e(1-e)\in (f_1,...,f_n)$.
\par Let $J=I\cap A$. By Lemma \ref{NQS} it is enough to find $j\in J^2$ and $F_i\in I_{1+j}$ such that $I_{1+j}=(F_1,...,F_n)$ with $f_i-F_i\in (I^2T)_{1+j}$. 
\par Let $B=\frac{A[T]}{(J^2[T],f_1)}$. Since $f_1$ is monic, $\dim(B)=\dim(\frac{A}{J^2})=\dim (\frac{A}{J})=\dim(\frac{A[T]}{I})\le n-1 $. Let `bar' denote modulo $(J^2[T],f_1)$. 
\par  In the ring $B$, we have $\ol I=(\ol f_2,...,\ol f_n)+\ol{I^2T}$. By Theorem \ref{IMMTT}, we get $\ol{h}_i\in \ol I$ such that $\ol I=(\ol h_2,...,\ol h_n)$ with $\ol f_i-\ol h_i\in \ol{I^2T}$ for all $i= 2,...,n$.
Hence we get, $I=(f_1,h_2,...,h_n)+J^2[T]$, where, $f_i-h_i\in (I^2T,J^2[T],f_1)$ for all $i= 2,...,n$. Note that by an elementary transformation we may further assume (we are not changing the notations $h_i$'s here) that $h_i-f_i\in I^2T+J^2[T]$ for all $i=2,...,n$. Define $F_i(T)=h_i(T)-h_i(0)+f_i(0)$. Then $F_i-f_i=(h_i(T)-f_i(T))-(h_i(0)-f_i(0))\in I^2T$ and $F_i\equiv h_i$ modulo $J^2[T]$ for all $i=2,...,n$. Therefore, we get $I=(f_1,F_2,...,F_n)+J^2[T]$ with $F_i-f_i\in I^2T$ for all $i=2,...,n$.
\par Again applying \cite{NMKL} we can find $s\in J^2A[T]$ with $s(1-s)\in (f_1,F_2,...,F_n)$. Let $I'=(f_1,F_2,...,F_n,1-s)$. Then we have $I\cap I'=(f_1,F_2,...,F_n)$ and $I'+J^2[T]=A[T]$. Since $I'$ contains a monic polynomial (namely $f_1$), then by \cite[Lemma 1.1, Chapter III]{Lam}  we can find $j\in J^2$ such that, $1+j\in I'\cap A$. We get $I_{1+j}=(f_1,F_2,...,F_n)_{1+j}$ with $f_i-F_i\in (I^2T)_{1+j}$. 
\smallskip

\paragraph{\textbf{Case - 2}}
In this case we shall assume that $n=2$. By \cite{SM} we may assume that $\dim(A[T]/I)+1=2$. Let $J=I\cap A$. By Lemma \ref{NQS} it is enough to find $j\in J^2$ and $h_i\in I_{1+j}$ such that $I_{1+j}=(h_1,h_2)$ with $f_i-h_i\in (I^2T)_{1+j}$.
\par  Since $I$ contains a monic polynomial, we have $\dim(A[T]/I)=\dim(A/J)=\dim(A/J^2)=1$. Let $C= \frac{A}{J^2}$,
then in the ring $C[T]$ we have $$\ol I=(\ol f_1,\ol f_2)+\ol{I^2T}.$$
Using Theorem \ref{IMMTT} we can find $g_i\in I$ such that $I=(g_1,g_2)+J^2[T]$, where $g_i-f_i\in I^2T+J^2[T].$ Let $h_i(T)=g_i(T)-g_i(0)+f_i(0)$. Then $h_i-g_i\in J^2[T]$ implies that $I=(h_1,h_2)+J^2[T]$ and $h_i-f_i\in I^2T.$ Now since $I$ contains a monic polynomial and $J$ is a proper ideal in $A$, the ideal $(h_1,h_2)$ contains a monic polynomial. By \cite{NMKL} there exists $s\in J^2A[T]$ with $s(1-s)\in (h_1,h_2)$. Let $I'=(h_1,h_2,1-s)$. Then $I'$ contains a monic polynomial with $I\cap I'=(h_1,h_2)$ and $I'+J^2[T]=A[T]$. By \cite[Lemma 1.1, Chapter III]{Lam}, there exists $j\in J^2$ such that $1+j\in I'\cap A$. Therefore we get $I_{1+j}=(h_1,h_2)_{1+j}[T]$ with $f_i-h_i\in (I^2T)_{1+j}$. This completes the proof.\qed\

\bc\label{QN}
Let $A$ be an affine algebra over  $\overline{\mathbb{F}}_p$. Let $I\subset A[T]$ be an ideal containing a monic polynomial such that $I(0)=(a_1,...,a_n)$. Assume that
$I=(f_1,...,f_n)+I^2$ with $f_i(0)-a_i\in I(0)^2$ for all $i$, where $n\ge \max\{(dim\frac{A[T]}{I}+1),2\}$. Then there exist $F_i\in I$ such that
$I=(F_1,...,F_n)$ with $F_i-f_i\in I^2$ and $F_i(0)=a_i$ for all $i$.
\ec

\proof Follows from \cite[Remark 3.9]{BR} and using Theorem \ref{1}.\qed 

\section{Improvement II: A relative version}
In \cite{MR} S. Mandal and R. Sridharan proved a relative version of Mandal's theorem quoted before. Their result has been crucial to the development of the Euler class theory. We now improve the bound of their result when the base ring is an affine algebras over  $\overline{\mathbb{F}}_p$.

\bt \label{Mtrick1}
Let $A$ be an affine algebra over  $\overline{\mathbb{F}}_p$. Let $I=I_1\cap I_2$, where $I_1, I_2\subset A[T]$ be two ideals such that:
\begin{enumerate}

	\item  $I_1$ contains a monic polynomial;
	\item  $I_1=(f_1,...,f_n)+I_1^2$, where  $n\ge \max\{(\dim(A[T]/I_1)+1),2\}$;
		\item $I_1+I_2=A[T]$;
	\item $I_2=I_2(0)A[T]$ (that is $I_2$ is extended from $A$).
\end{enumerate}
Suppose that there exist $a_i\in I(0)$ with $a_i-f_i(0)\in I_1(0)^2$ such that $$I(0)=(a_1,...,a_n)$$ for all $i$. Then there exist $h_i\in I$ with $h_i(0)=a_i$ such that $$I=(h_1(T),...,h_n(T)).$$
\et
\proof Let $J_1=I_1\cap A$. Since $I_1$ contains a monic polynomial and $I_2$ is extended from $A$, by \cite[Lemma 1.1, Chapter III]{Lam} we can find $s\in J_1$ and $t\in I_2(0)$ such that $s+t=1$. Note that in the ring $A_t[T]$, we have $I_1A_t[T]=(f_1,...,f_n)+I_1^2A_t[T]$ and $I_1(0)A_t=I(0)A_t=(a_1,...,a_n)A_t$ with $f_i(0)-a_i\in I_1(0)^2A_t $. Therefore, by Corollary \ref{QN} there exist $g_i\in I_1A_t[T]$, for $i=1,...,n$ such that $I_1A_t[T]=(g_1,...,g_n)$ with $g_i(0)=a_i$. Now consider the following exact sequences

$$\begin{tikzcd}
0 \arrow{r} & (K_1)_s \arrow{r} & A^n_{st}[T] \arrow{rr}{(g_1,...,g_n)} && IA_{st}[T](=I_1A_{st}[T]=A_{st}[T]) \arrow{r} & 0\\ 
0 \arrow{r} & (K_2)_{st} \arrow{r} & A^n_{st}[T] \arrow{rr}{(a_1,...,a_n)\otimes A_{st}[T]} & & IA_{st}[T](=I(0)_{st}A[T]=A_{st}[T]) \arrow{r} & 0. 
\end{tikzcd}$$ 
where $K_1$ is the kernel of the map from $A^n_t[T]\to I_1A_{t}[T]$ induced by $(g_1,...,g_n)$. Moreover, it follows from the proof of Theorem \ref{1} that $g_1$ can be taken as a monic polynomial. Hence by \cite{Q} the module $K_1$ is extended from $A_t$. The module $K_2$ is the kernel of the map from $A^n[T]\to I(0)A[T] $ induced by $(a_1,...,a_n)\otimes A[T]$, which is also extended from $A$, as the map itself is extended. Then using \cite[Lemma 2]{P} we can find $\alpha(T) \in SL_n(A_{st}[T])$ such that $\alpha(0)=\text{Id}$ and $(g_1,...,g_n)\alpha(T)=(a_1,...,a_n)$. Therefore, using Quillen's Splitting Lemma \cite{Q} we get $\alpha(T)=(\alpha_1(T))_s(\alpha_2(T))_t$, where $\alpha_1(T)\in \GL_n(A_{s}[T]) $ and  $\alpha_2(T)\in\GL_n(A_{t}[T]) $. Then a standard patching argument completes the proof.\qed

\section{Applications I: Addition and Subtraction Principles}
\bp\label{apA}(Addition Principle)
Let $A$ be an affine algebra over  $\overline{\mathbb{F}}_p$ of dimension $d$. Suppose that $I_1,I_2\subset A$ be two co-maximal ideals of height $n$, where $2n\ge d+2$. Let $I_1=(a_1,...,a_n)$, $I_2=(b_1,...,b_n)$ and $I=I_1\cap I_2$. Then there exist $c_i\in I$ such that $I=(c_1,...,c_n)$ with $c_i-a_i\in I_1^2$ and $c_i-b_i\in I_2^2$. 
\ep
\proof Without loss of generality we may assume that $d>n$, as the case $d<n$, follows from \cite[Corollary 3]{NMK} and $d=n$, follows from \cite{SMBB3}. Note that we can always perform elementary transformations on $(a_1,...,a_n)$ and $(b_1,...,b_n)$. Let $B=A/(b_1,...,b_n)$ and 'bar' denote going modulo $(b_1,...,b_n)$. Note that $(\ol{a}_1,...,\ol{a}_n)\in \Um_n(B)$. Since $n\ge d-n+2\ge\dim(B)+2$, we shall have $(\ol a_1,...,\ol a_n)\sim_{E_n(R)}(\ol 1,...,\ol 0)$. Adding suitable multiples of $a_n$ to $a_i$'s, by a theorem of Eisenbud-Evans \cite{EE} we may further assume that $\hh(a_1,...,a_{n-1})=n-1$. Since $\ol a_n=0$, we still have $(\ol a_1,...,\ol a_n)\sim_{E_n(R)}(\ol 1,...,\ol 0)$. Therefore, we get $(a_1,...,a_{n-1})+I_2=A$.
\par Let $C=A/(a_1,...,a_{n-1})$ and `tilde' denote going modulo $(a_1,...,a_{n-1})$. Since  $(a_1,...,a_{n-1})+I_2=A$, we have $( \wt b_1,..,\wt b_n)\in \Um_n(C) $. Also note that $n\ge (d-n+1)+1\ge \dim(C)+1$. Therefore, by \cite[Corollary 17.3]{SV} we get $(\wt b_1,..,\wt b_n)\sim_{E_n(R)}( \wt 1,..., \wt 0)$. We may further assume that  $( \wt b_1,..,\wt b_n)=( \wt 1,..., \wt 0)$. Again, as before, without altering the assumption $ \wt b_n=\wt 0$ we may also assume that $\hh(b_1,...,b_{n-1})=n-1$. Hence we get:
\\(i) $(a_1,...,a_{n-1})+(b_1,...,b_{n-1})=A$,
\\(ii) $\hh(a_1,...,a_{n-1})=\hh(b_1,...,b_{n-1}) =n-1$.
\par Now define $J_1=(a_1,...,a_{n-1},a_n+T)\subset A[T] $,  $J_2=(b_1,...,b_{n-1},b_n+T)\subset A[T]$ and let $J=J_1\cap J_2$ (note that $J$ contains a monic polynomial). Since $J_1+J_2=A[T]$, by the Chinese Remainder Theorem $J/J^2\cong J_1/J_1^2\oplus J_2/J_2^2$. Therefore, we can find $g_i\in J$, for $i=1,...,n$ such that $J=(g_1,...,g_n)+J^2$ with $g_i-a_i\in J_1^2$, $g_i-b_i\in J_2^2$, for $i=1,...,n-1$, $g_n-a_n-T\in J_1^2$ and $g_n-b_n-T\in J_2^2$. Also $A[T]/J\cong A[T]/J_1\oplus A[T]/J_2$ implies that
$\dim(A[T]/J)=\max \{\dim(A[T]/J_1),\dim(A[T]/J_2)\}=\max \{\dim(A/(a_1,...,a_{n-1})),\dim(A/(b_1,...,b_{n-1}))\}=d-n+1\le n-1$.
\par Therefore, by \cite[Theorem 3.2]{MKD} we can find $h_i \in J$, $i=1,...,n$ such that $J=(h_1,...,h_n)$ with $h_i-g_i\in J^2$. We define $c_i$ to be the evaluation of $h_i(T)$ at $T=0$, that is, $c_i:=h_i(0)$ for all $i=1,...,n$. Then 
$I_1\cap I_2=(c_1,...,c_n)$ with $c_i-a_i\in I_1^2$ and $c_i-b_i\in I_2^2$, for $i=1,...,n$. \qed

\bp\label{SPSFM}(Subtraction Principle)
Let $A$ be an affine algebra over  $\overline{\mathbb{F}}_p$ of $\dim(A)=d$. Suppose that $I_1,I_2\subset A$ be two comaximal ideals of height $n$, where $2n\ge d+2$. Let $I_1=(a_1,...,a_n)$ and $I=I_1\cap I_2=(c_1,...,c_n)$ with $c_i-a_i\in I_1^2$ for all $i=1...,n$. Then there exist $b_i\in I$ such that $I_1=(b_1,...,b_n)$ with $c_i-b_i\in I_2^2$.
\ep
\proof Using the same arguments as used in the proof of Proposition \ref{apA} we may establish the followings:

\begin{enumerate}
	\item $\hh(a_1,...,a_{n-1})=n-1$;
	\item $a_1-1\in I_2$ and $a_i\in I_2$ for all $i>1$;
	\item Replacing $a_n$ by $a_1+a_n$, we may further assume that $a_n-1\in I_2$.
	
	\end{enumerate}

\par Define $J_1=(a_1,...,a_{n-1},a_n+T)$, $J_2=I_2A[T]$ and $J=J_1\cap J_2$. Then we get:
\begin{enumerate}
	\item $J_1+J_2=A[T]$;
	\item  $J_1$ contains a monic polynomial;
	\item  $J_1=(a_1,...,a_{n-1},a_n+T)+J_1^2$ with $\dim(A[T]/J_1)+1=\dim(A/(a_1,...,a_{n-1}))+1\le d-n+2\le n$ ;
	\item  $J_2$ is an extended ideal;
	\item $J(0)=I_1\cap I_2=I= (c_1,...,c_n)$ with $c_i-a_i\in I_1^2=J_1(0)^2$ for all $i=1...,n$.
	
\end{enumerate}

\par Therefore, applying Theorem \ref{Mtrick1}, we can obtain $J=(h_1,...,h_n)$ such that evaluating $h_i(T)$ at $T=0$ match with $c_i$ for all $i=1...,n$. We define $b_i$ to be the evaluation of $h_i(T)$ at $T=1-a_n$, that is, $b_i:=h_i(1-a_n)$ for all $i=1,...,n$. Then we get $I_2= J_2(1-a_n)=J_1(1-a_n)\cap J_2(1-a_n)=J(1-a_n)=(b_1,...,b_n)$ with $c_i-b_i=h_i(0)-b_i=h_i(1-a_n)-b_i=0$ modulo $I_2^2$. This completes the proof.  \qed


\section{Application II: The Euler class group and the Euler class}\label{7}

Let $A$ be an affine algebra over $\overline{\mathbb{F}}_p$ such that $\dim(A)=d.$ Henceforth we shall assume that $d\ge 3$ and $n$ is an integer satisfying $2n\ge d+2$.
\par Let $J$ be an ideal of height $n$ such that $J/J^2$ is generated by $n$ elements. Let $\alpha$ and $\beta$ be two surjections from $(A/J)^n$ to $J/J^2.$ We say that $\alpha$ and $\beta$ are related if there exists an elementary automorphism $\sigma $ of $(A/J)^n$ such that $\alpha\sigma=\beta$. This is an equivalence relation on the set of surjections from $(A/J)^n$ to $J/J^2$. Let $[\alpha]$ denote the equivalence class of $\alpha$. If $\ol{a}_1,...,\ol{a}_n$ generate $J/J^2$, we obtain a surjection $\alpha:(A/J)^n\surj/J^2$, sending $\ol{e}_i $ to $\ol{a}_i$. We say $[\alpha]$ is given by the set of generators $\ol{a}_1,...,\ol{a}_n$ of $J/J^2$.
\medskip

\bd  \label{DE}\par Let $G$ be the free abelian group on the set $B$ of pairs $(J,\omega_J)$, where:\\
$(i)$ $J\subset A$ is an ideal of height $n$,\\
$(ii)$ $\Spec(A/J)$ is connected,\\ 
$(iii)$ $J/J^2$ is generated by $n$ elements and \\
$(iv)$ $\omega_J:(A/J)^n\surj J/J^2$ is an equivalence class of surjections $\alpha:(A/J)^n\surj J/J^2$.
\medskip

\par Let $J\subset A$ be a proper ideal. Then $J=J_1\cap J_2\cap ...\cap J_r$, where $J_i$'s are proper, pairwise co-maximal and $\Spec(A/J_i)$ is connected for all $i=1,...,r$. It was proved in \cite[Lemma 4.1]{BR3} that such a decomposition is unique. We shall say that $J_i$ are the connected components of $J$.
\par Let $J\subset A$ be an ideal of height $n$ such that $J/J^2$ is generated by $n$ elements. Let $J=\cap J_i$ be the decomposition of $J$ into its connected components. Then note that  for every $i$, $\hh(J_i)=n$ and by the Chinese Remainder Theorem  $J_i/J_i^2$ is generated by $n$ elements. Let $\omega_J:(A/J)^n\surj J/J^2$ be a surjection. Then in a natural way $\omega_J$ gives rise to surjections $\omega_{J_{i}}:(A/J_i)^n\surj J/{J_i}^2$. We associate to the pair $(J,\omega_J)$, the element $\sum (J_i,\omega_{J_i})$ of $G$.
\par Let $H$ be the subgroup of $G$ generated by the set $S$ of pairs $(J,\omega_J)$, where  $\omega_J:(A/J)^n\surj J/J^2$ has a surjective lift and denote it by $\theta:A^n\surj J $. Then we define the quotient group $G/H$ as the $n$-th Euler class group of $A$ and denote it by ${E^n(A)}$.
\ed
\bt\label{DEG}
Let $A$ be an affine algebra over  $\overline{\mathbb{F}}_p$ of dimension $d$ and $n$ be an integer satisfying $2n\ge d+2$. Let $I\subset A$ be an ideal of height $n$ be such that $I/I^2$ is generated by $n$ elements and $\omega_I:(A/I)^n\surj I/I^2 $ be an equivalence class of surjections. Suppose that the image of $(I,\omega_I)$ is zero in the Euler class group $E^n(A)$ of $A$. Then $I$ is generated by $n$ elements and $\omega_I$ can be lifted to a surjection $\theta:A^n\surj I.$
\et

To prove Theorem \ref{DEG}, we shall need the following lemma whose proof can be found in \cite[Lemma 4.1]{MPHIL}.

\bl\label{BajeLemma}
Let $F$ be a
free abelian group with basis $\{e_i\}_{i\in\MI} $ and $\sim$ be an equivalence relation on $\{e_i\}_{i\in \MI}$. Define $x \in F$ to
be ``reduced'' if $x = e_1 + . . . + e_r$ and $e_i\not= e_j$ ($i \not= j$). For $x\in F$ with $x = e_1 + . . . + e_r$ we define ``support" of $x$ to be the set $supp(x) = \{e_1,...,e_r\}$. Define, $x \in F$ to be ``nicely
reduced' if $x = e_1 + . . . + e_r$ and $e_i\not= e_j$ ($i \not= j$) and such that no $e_i$ belongs to the equivalence class of other $e_j$ for $i,j=1,...,r$ and $i \not= j$. Let $S \subset F$
be such that:
\begin{enumerate}

\item Every element of $S$ is nicely reduced;
\item  Let $x,y\in F$ such that $x, y$ and $ x + y$ are nicely reduced. Then if any two of $x,y,x+y$ are in $S$, then so is the third one;
\item  $x \in F$ and $x \not\in  S$ be nicely reduced. Let $K \subset \MI$ be a finite set. Then there exists $y \in  F$ such that:
\begin{enumerate}[(i)]

\item $y$ is nicely reduced;
\item $x + y \in S$;
 \item  $y + e_k$ is nicely reduced for all $k \in K$.
\end{enumerate}

\end{enumerate}
Let $ H$ be a subgroup of $F$ generated by the set $S$ . Then if $x \in H$ is nicely reduced then $x \in S$.
\el

\noindent
Proof of  Theorem \ref{DEG}.  We take $F$ to be the free abelian group generated by the set $B$, as defined in \ref{DE}. Define
a relation `$\sim$' on B as $(J,\omega_J)\sim (I,\omega_I)$ if $I=J$. Then it is an equivalence relation.
\\
Let $S\subset G$ be as in \ref{DE}. In view of the above lemma, it is enough to show that the
three conditions in Lemma \ref{BajeLemma} are satisfied. Condition $(i)$ is clear, almost from the definition.
The addition and subtraction principles (\ref{apA} and \ref{SPSFM}) will yield condition $(ii)$.  Finally,
applying the moving lemma \cite[Corollary 2.4]{BR3}, it is clear that $(iii)$ is also satisfied. \qed

\subsection{The Euler class of stably free modules}\label{8}
Recall that \emph{`1-stably free modules'} are those which are given by unimodular rows. Let $A$ be a smooth affine algebra over  $\overline{\mathbb{F}}_p$ of dimension $d\ge 3$. In this section we shall assign an $1$-stably free module $P$ of rank $n$ to its \emph{`$n$-th Euler class'} and show that $P$ has a unimodular element if and only if its Euler class is trivial. Let $v=(v_0,...,v_n)\in \Um_{n+1}(A)$, where $2n\ge d+2$. Note that the case $d=n$ is done in \cite{SMBB3}. So without loss of generality we may further assume that $d\ge n+1$. Let $\displaystyle P=A^{n+1}/(\sum_{0}^{n} v_i e_i)$, where $e_0,\cdots,e_n$ is the canonical basis  of $A^{n+1}$.    Therefore, we have $P\oplus A\cong A^{n+1}$. Let $p_i$ denote the image of $e_i$ in $P$, for $i=0,...,n$ with $p=(p_0,...,p_n)$. Then $P=\sum_{i=0}^{n}Ap_i$ and $\sum_{i=0}^{n}v_ip_i=0$. We shall define a map $\Um_{n+1}(A)\to E^n(A)$ for all $n$, which satisfies $2n\ge d+2$ and shall assign $P$ with an element of the group $E^n(A)$ in the following way:
\par Let $\lambda:P\surj J$ be a surjection, where $J\subset A$ is an ideal of height $n.$ Since $P$ is a $1-$stably free $A$-module of rank $n\ge d-n+2\ge \dim(A/J)+2 $, $P/JP$ is a free $A/J$-module. Note that $\lambda:P\surj J$ will induce a  surjection $\ol{\lambda}:P/JP\surj J/J^2$ and since $P/JP$ is free, $J/J^2$ is generated by $n$ elements. 
\par Let `bar' denote going modulo $J$. Since $\dim(A/J)\le d-n$, we have $n+1\ge d-n+1 \ge \dim(A/J)+1\ge \sr(A/J)+1$. Therefore, by \cite[Corollary 17.3]{SV} $\ol v\sim_{E_{n+1}(A)} \ol e_1$.  Hence there exists $\epsilon \in E_{n+1}(A)$ with  $\ol \epsilon\ol e_1^T=\ol v^T$ and therefore, we have $\ol p \ol \epsilon =(\ol 0, \ol u_1,...,\ol u_n)$. This implies $\{\ol u_1,...,\ol u_n\}$ forms a basis of the free $A/J$-module $P/JP$. Let $\omega_J$ be the surjection induced by $\lambda$, sending each $\ol u_i \to \ol {\lambda (u_i)}$, for $i=1,...,n$. 
\par We assign $e(P,v,p):= (J,\omega_J)\in E^n(A)$. We check below that this assignment is well defined.

\bt \label{w1} With the same notations as above suppose that there exists another surjection $\lambda':P\surj J'$, where $J'\subset A$  be an ideal of height $n$ and we obtain $\omega_{J'}$ in the same way as discuss earlier. Then $(J,\omega_J)=(J',\omega_{J'})$ in $E^n(A)$.
\et
\proof By \cite[Lemma 5.1]{BR3} there exists an ideal $I\subset A[T]$ of height $n$ and a surjection $\lambda(T):P[T]\surj I $ such that $\lambda(0)=\lambda$ and $\lambda(1)=\lambda'$. Let $N=(I\cap A)^2$. Then $n-1\le \hh(N)\le n$ and by our assumption $\dim(A/N)\ge 2$. Therefore, by  \cite[Corollary 17.3]{SV} we have $\sr(A/N)\le \dim(A/N)\le d-n+1<d-n+3\le n+1. $ Hence there exists $\alpha \in E_{n+1}(A)$ with $\ol\alpha\ol e_1^t=\ol v^t$ and we have $\ol p \ol \alpha =(\ol 0, \ol w_1,...,\ol w_n)$. Therefore,$\{\ol w_1,...,\ol w_n\}$ forms a basis of the free $A/N$-module $P/NP$ and also a basis of the free $A[T]/I$-module $P[T]/IP[T]$. Therefore, as earlier we obtain $\{\ol \lambda(T)(\ol w_1),...,\ol \lambda(T)(\ol w_n)\}$ as a set of generators of $I/I^2$ and setting $T=0$ and $1$ we obtain generators of $J/J^2$ and $J'/J'^2$, receptively and hence the corresponding surjections $\omega_J:(A/J)^n\surj J/J^2$ and $\omega_J':(A/J')^n\surj J'/J'^2$. 
\par Then by \cite[Proposition 5.2]{BR3} there exists an ideal $K\subset A$ of height $n$, comaximal with $J$ and $J'$ and a surjection $\omega_K:(A/K)^n\surj K/K^2$ such that $(J,\omega_J)+(K,\omega_K)=(J',\omega_J')+(K,\omega_K)$ in $E^n(A)$. As a consequence $(J,\omega_J)=(J',\omega_J')$.          \qed

\bt \label{w2} Suppose that there exists $\epsilon' \in E_{n+1}(A)$ with  $\ol \epsilon'\ol e_1^T=\ol v^T$ such that $\ol p \ol \epsilon' =(\ol 0, \ol u'_1,...,\ol u'_n)$. Let $\omega'_J$ be a surjection induced by $\lambda$, sending each $\ol u'_i \to \ol {\lambda (u'_i)}$, for $i=1,...,n$. Then  $(J,\omega_J)=(J,\omega'_J)$. 
\et

Theorem \ref{w2} follows from the following lemma.
\bl\label{L3}
Suppose that there exists $ \epsilon'\in E_{n+1}(A)$ with $\ol \epsilon'\ol e_1^T=\ol v^T$ such that  $\ol p \ol \epsilon' =(\ol 0, \ol u'_1,...,\ol u'_n)$. Then there exists $\theta \in E_n(A/J)$ such that $(\ol u_1,...,\ol u_n)\theta =(\ol u'_1,...,\ol u'_n)$.
\el
\proof Since  $\ol \epsilon\ol e_1^T=\ol v^T=  \ol \epsilon'\ol e_1^T$, we have  $\ol \epsilon^{-1}\epsilon'\ol e_1^T=\ol e_1^T$ and $\epsilon^{-1}\epsilon'\in E_{n+1}(A)$. Therefore, there exists $\theta \in SL_n(A/J)\cap E_{n+1}(A/J)$ such that $(\ol u_1,...,\ol u_n)\theta =(\ol u'_1,...,\ol u'_n)$. Since $n\ge d-n+2>\dim(A/J)=\sr(A/J)$, by \cite[Theorem 3.2]{V}  $\theta \in E_n(A/J)$.\qed
\smallskip 
\par Hence the assignment of $(J,\omega_J)$ associated to $(P,v,p)$ is well defined and we shall denote it by $e(v)\in E^n(A)$. From the definition of $E^n(A)$ it follows that, for any two unimodular rows $u,v\in \Um_{n+1}(A)$, if $u\sim_{E_{n+1}(A)}  v$, then $e(v)=e(u)$. Therefore, we obtain a map $$e:\Um_{n+1}(A)/E_{n+1}(A)\to E^n(A).$$

\bt \label{NS}
Suppose that $v$ and $P$ be as defined before. Then $P$ has a unimodular element if and only if $e(v)=0$ in $E^n(A)$.
\et
\proof First we assume that $P$ has a unimodular element, say $u\in \Um(P)$. Therefore, we have $P\cong Q\oplus A{u}$. Let $\lambda:P\surj J$ be a surjection, where $J\subset A$ be an ideal of height $n$. By (\cite{EE} also see, \cite[Corollary 2.13]{SMBB3}) without loss of generality we may assume that, $\hh(\lambda(Q))=n-1.$ Let $N=\lambda(Q)$ and `bar' denote going modulo $N$. Since $n+1\ge d-n+3= \dim(A/N)+2$, there exists $\epsilon\in E_{n+1}(A)$ with $\ol\epsilon\ol e_1^t=\ol v^t$. Therefore, we have $\ol p \ol \epsilon =(\ol 0, \ol w_1,...,\ol w_n)$. Hence $\{\ol w_1,...,\ol w_n\}$ will form a basis of the free $A/N$-module $P/NP$.
\par Let $\ol u=\sum_{i=1}^{n}\ol a_i\ol w_i$. Since $\ol u$ is a unimodular element of the free $A/N$-module $P/NP$, we have $(\ol a_1,...,\ol a_n)\in \Um_n(A/N)$. Since $\dim(A/N)=d-n+1\ge 2$ we have, $n\ge \dim(A/N)+1\ge \sr(A/N)+1$. This implies that $(\ol a_1,...,\ol a_n)\sim_{E_{n+1}(A/N)} \ol e_1$. Hence $(\ol w_1,...,\ol w_n)$ can be taken by an elementary automorphism to a basis $\{\ol u,\ol t_2,...,\ol t_n\}$, where $\ol t_i\in Q/NQ$, for $i=2,...,n$. Let $\ol \lambda(t_i)=\ol b_i$, for $i=2,...,n$. Therefore, we get $J=(b_2,...,b_n,c)+J^2$ and $N=(b_2,...,b_n)+N^2$, where $c=\lambda(u)$. Then by \cite{NMKL} there exists $e\in N^2$ such that $e(1-e)\in (b_2,...,b_n)$. Let $N=(b_2,...,b_n,e)$ and $J=(N,c)=(b_2,...,b_n,e+(1-e)c)$. Therefore, we obtain a lift of $J=(b_2,...,b_n,c)+J^2$, which implies that $e(v)=0$ in $E^n(A).$
\par Conversely, we assume that $e(v)=0$ in $E^n(A)$. Let $\lambda:P\surj J$, be a surjection, where $J\subset A $ be an ideal of height $n$. Let `bar' denote going modulo $J$. Since $n+1\ge \dim(A/J)+2$, there exists $\alpha \in E_{n+1}(A)$ such that $\ol v^t=\ol \alpha\ol e_1^t$. Therefore, we have $\ol p\ol \alpha= (\ol 0,\ol u_1,...,\ol u_n)$ and a basis $\{\ol u_1,...,\ol u_n\}$ of the free $A/J$-module $P/JP$. Hence we get $J=( {\lambda(u_1)},...,{\lambda(u_n)})+J^2$. Since $e(v)=0$, we can obtain $J=(b_1,...,b_n)$ with $\lambda(u_i)-b_i\in J^2$, for $i=1,...,n$. Moreover using a theorem of Eisenbud-Evans \cite{EE} we can find $c_1,...,c_{n-1}\in A$ such that $\hh(b_1+c_1b_n,...,b_{n-1}+c_{n-1}b_n)=n-1$. Let $d_i=b_i+c_ib_n$, for $i=1,...,{n-1}$ and $d_n=b_n$.
\par Let $I=(d_1,...,d_{n-1},d_n+T)$ and `tilde' denote going modulo $I$. Then we have $P[T]=(A[T])^{n+1}/(v)$. We see that $\dim(A[T]/I)=\dim(A/(d_1,...,d_{n-1}))=d-n+1\le n-1<n+1$ and hence we can find $\Gamma(T)\in E_{n+1}(A[T])$ such that $\widetilde{\Gamma(T)} \widetilde{e_1}^t= \widetilde{v}^t$. Therefore, we get $\widetilde{p}\widetilde{\Gamma(T)}=(\widetilde{0},\widetilde{u_1(T)},...,\widetilde{u_n(T)})$ and $\{\widetilde{u_1(T)},...,\widetilde{u_n(T)}\}$ will form a basis of the free $A[T]/I$-module $P[T]/IP[T]$. Evaluating at $T=0$, we get $\ol p\ol \Gamma(0)=(\ol 0,\ol u_1(0),...,\ol u_n(0))$ with $\Gamma(0)\in E_{n+1}(A)$. 
\par Therefore, by Lemma \ref{L3} we can find $\theta_1 \in E_n(A/J)$ such that, $(\ol u_1,...,\ol u_n)= (\ol u_1(0),...,\ol u_n(0))\theta_1$. Hence there exists $\theta_2\in E_{n+1}(A/J)$ such that, $$(\ol{u_1+c_1u_n},...,\ol{u_{n-1}+c_{n-1}u_n},\ol{u_n})=(\ol{u_1(0)},...,\ol{u_n(0)})\theta_2$$
Note that $\ol{\lambda(u_i+c_iu_n)}=\ol{d_i}$ for all $i=1,...,n-1$ and $\ol{\lambda(u_n)}=\ol{d_n}$. Also note that $A/J\cong (A[T]/I)/(t)$, where $t=\widetilde{T}$. Therefore, the map $A[T]/I\surj A/J$ is surjective and so is $E_{n}(A[T]/I)\surj E_{n}(A/J)$. Hence we can find $\tau(T) \in E_n(A[T])$, which is a lift of $\theta_2$. Let 
$$(\widetilde{u_1(T)},...,\widetilde{u_n(T)})\widetilde{\tau(T)}=(\widetilde{w_1(T)},...,\widetilde{w_n(T)}).$$
Since $(\widetilde{u_1(T)},...,\widetilde{u_n(T)})$ is a basis of the free $A[T]/I$-module $P[T]/IP[T]$, $\{\widetilde{w_1(T)},...,\widetilde{w_n(T)}\}$ is also a basis. Define a surjection $\theta:P[T]/IP[T]\surj I/I^2 $, sending $\widetilde{w_i(T)}\to {d_i}$ for $i=1,...,n-1$ and $\widetilde{w_n(T)}\to {d_n+T}$. Then $\widetilde{\tau(0)}=\theta_2$, implies that $\theta(0)=\ol \lambda$. Since we have $d\ge n+2$ this gives us  $\dim(A[T]/I)+1\le d-n+2\le n$, by Corollary \ref{QN} there exists a surjective lift $\Theta$ of $\theta :P[T]\surj  I$ which matches at $T=0$. Since $T+d_n\in I$, setting $T=1-d_n$, we obtain a surjection from $\gamma=\Theta(1-d_n): P\surj A$, which completes the proof.\qed

\bl\label{VS}
Let $n$ be an even integer. Let $J=(a_1,...,a_n)$ be an ideal of height $n$ and $u$ be a unit modulo $J$. Let $\omega_J:(A/J)^n\surj J/J^2$ be given by the set of generators $\ol{ua_1},\ol{a_2},...,\ol{a_n}$ of $J/J^2$. Let $v\in A$ be such that $1-uv\in J$. Then $e(v,a_1,...,a_n)=(J,\omega_J)$ (Note that $(v,a_1,...,a_n)\in \Um_{n+1}(A)$).
\el
\proof Let $Q=A^{n+1}/(v,a_1,...,a_n)$ and $q_i$ denote the image of $e_i$ in $P$, for $i=0,...,n$. Let $\mu:Q\surj J$ be a surjection sending $q_0\to 0$ and $q_i\to a_{i+1}$ if $i$ is odd and $q_i\to -a_{i-1}$ if $i$ is even. Therefore, modulo $J$ we get $(v,a_1,...,a_n)=(v,0,...,0)$. By Whitehead's Lemma, the diagonal matrix given by $\diag(\ol v, \ol u, \ol 1,\ol 1...,\ol 1)\in E_n(A/J)$. Hence $e(v,a_1,...,a_n)=(J,\omega_J)$, where $\omega_J$ is given by the set of generators $\{\ol{ua_2},-\ol a_1, \ol a_4,-\ol a_3,...,$ $\ol a_n,-\ol a_{n-1}\}$ of $J/J^2$. Applying Whitehead's Lemma again, we see that $\omega_J$ is given by the set of generators $\ol{ua_1},\ol a_2,...,\ol a_n$ of $J/J^2$.  \qed

\smallskip

Recall that the stable dimension of a ring $R$ is defined by $sdim(R)=\sr(R)-1$. The next theorem can be found in \cite[Theorem 4.1]{VdKM}.
\bt \label{TWVK} Let $R$ be a commutative ring, $n\ge 3$ and $sdim(A)\le 2n-4$. Then the universal weak Mennicke symbol $wms:\Um_n(A)/E_n(A)\to WMS_n(R)$ is bijective. Therefore, $\Um_n(A)/E_n(A)$ has a group structure. 
\et 
 Since  $A$ is an affine algebra over $\overline{\mathbb{F}}_p$ of dimension $d\ge 3$, by \cite[Corollary 17.3]{SV} we have $sdim(A)\le d-1$. Therefore, by Theorem \ref{TWVK} for all $n$ satisfying $2n\ge d+1$, we have a group structure on $\Um_{n+1}(A)/E_{n+1}(A)$.
\bt
For all $n$ satisfying $2n\ge d+2$, the map $e:\Um_{n+1}(A)/E_{n+1}(A)\to E^n(A)$ is a group homomorphism.
\et
\proof If $n$ is odd, since every $1-$stably free module of odd rank has a unimodular element, using Theorem \ref{NS} it follows that $e$ is the zero map. In this case there are nothing to prove. So shall assume that $n$ is even. Without loss of generality we may further assume that $A$ is a domain. Since $2(n+1)\ge d+4\ge d+3$, by \cite[Lemma 3.2]{MNLvdK} it is enough to prove that if $(x,a_1,...,a_n)$ and $(1-x,a_1,...,a_n)$ are unimodular rows then we have
$$e(x,a_1,...,a_n)+e(1-x,a_1,...,a_n)=e(x(1-x),a_1,...,a_n).$$
Let $y=1-x$. We may assume that $xy\not= 0$. Let `bar' denote going modulo $xy$. Then adding a suitable multiple of $\ol a_1$ to $\ol a_i$ for $i=2,...,n$ we may assume that $\hh(\ol a_2,...,\ol a_n)=n-1$ and hence $\hh(x,a_2,...,a_n)=n.$ Therefore, we may assume $\hh(y,a_2,...,a_n)=n=\hh(x,a_2,...,a_n).$
\par Let $b_1\in A$ be such that $1+a_1b_1\in (xy,a_2,...,a_n)$. Now since $(x,a_1,a_2...,a_n)\sim_E(-a_1,x,a_2...,a_n) $, we have $e(x,a_1,a_2...,a_n)=e(-a_1,x,a_2...,a_n) $.
\par Since $n$ is even by Lemma \ref{VS}, $e(-a_1,x,a_2...,a_n)=e(x,a_1,a_2...,a_n)=(J_1,\omega_{J_1})$, where $J_1=(x,a_2,...,a_n)$ and $\omega_{J_1}$ is given by the set of generators $\ol{b_1x},\ol a_2,...,\ol a_n$ of $J_1/J_1^2$. Similarly we get $e(y,a_1,a_2...,a_n)=(J_2,\omega_{J_2})$, where $J_2=(y,a_2,...,a_n)$ and $\omega_{J_2}$ is given by the set of generators $\ol{b_1y},\ol a_2,...,\ol a_n$ of $J_2/J_2^2$ and $e(xy,a_1,a_2...,a_n)=(J_3,\omega_{J_3})$, where $J_3=(xy,a_2,...,a_n)$ and $\omega_{J_3}$ is given by the set of generators $\ol{b_1xy},\ol a_2,...,\ol a_n$ of $J_3/J_3^2$. Since $J_3=J_1\cap J_2$ and $x+y=1$ we have $(J_3,\omega_{J_3})=(J_1,\omega_{J_1})+(J_2,\omega_{J_2})$. This completes the proof.\qed

\section{Application III: On a question of Roitman}

\bt\label{QORK}
Let $R$ be an affine $\k-$algebra of dimension $d\geq 2$ and $P$ be a $1$-stably free $R[T]$-module of rank $d-1$. Assume that $P\otimes R(T)$ has a unimodular element. Then $P$ also has a unimodular element.
\et

\proof
If $d=2$, then the result follows trivially. If $d=3$, then since $P$ is of rank two with trivial determinant, we observe that $P\ot R(T)$ is actually free. Then, by the Affine Horrocks Theorem, $P$ is free. Therefore, we assume that $d\geq 4$. As we shall apply the Euler class theory developed in Section \ref{8}, we take $2(d-1)\geq (d+1)+2$, implying that $d\geq 5$. Apparently, the only case that seems to be left out is $d=4$ and $\rank(P)=3$. But it is easy to see that $1$-stably free modules of odd rank always have a unimodular element. Therefore, all the cases will be covered by this theorem once we complete the following arguments (with $d\geq 5$).

Since $R$ contains an infinite field, namely, $\k$, we can follow the arguments of \cite[Lemma 3.1]{BRS01} and obtain an $R[T]$-linear surjection
$\lambda:P\sur I$, where $I\subset R[T]$ is an ideal of height $d-1$ and $I$ contains a monic polynomial. Let $P$ correspond to the unimodular row $v\in \Um_{d}(R[T])$. Using $\lambda$ we can compute the Euler class of $P$ (or $v$), as in Section \ref{8} and obtain
$$e(v)=(I,\omega_I)\in E^{d-1}(R[T]).$$

As $I$ contains a monic polynomial, it follows from \cite[Theorem 3.2]{MKD} that $\omega_I:(R[T]/I)^{d-1}\sur I/I^2$ has a surjective lift $\theta:R[T]^{d-1}\sur I$. In other words, $(I,\omega_I)=0$ in $E^{d-1}(R[T])$. Consequently, $e(v)=0$ and by Theorem \ref{NS}, $P$ has a unimodular element.
\qed

\smallskip

\rmk One can prove the above theorem for arbitrary stably free modules of rank $d-1$ with the help of the arguments given in \cite[Section 9]{STMKK}. It is worth noting that in \cite{STMKK}, the assignment of an `Euler cycle' to a stably free module $P$ does not require any hypothesis related to the `smoothness' of the ring. Additionally, it follows from Theorem \ref{NQ} that the vanishing of any Euler cycle assigned to $P$ is sufficient for $P$ to split off a free summand of rank one.

\appendix

\section{A ``projective" version of Theorem \ref{1}}

Theorem \ref{1} was proved with the free module $A^n$. Here, in place of $A^n$ we take a projective $A$-module $P$ of rank $n$.

\bt\label{NQ}
Let $A$ be an affine algebra over $\ol{\mathbb{F}}_p$ and $I\subset A[T]$ be any proper ideal containing a monic polynomial. Suppose that $P$ is a projective $A$-module of rank $n$, where $n\ge \max\{ (\dim A[T]/I+1),2\}$.
Then any surjective map $\ol{\phi}:P[T]\surj I/I^2T$  lifts  to a surjective map $\phi:P[T]\surj I$.
\et

\proof Let $J=I\cap A$. Since $I$ contains a monic polynomial $\dim(A[T]/I)=\dim(A/J)\le n-1.$ Since $\rank(P/JP)=n> \dim(A/J)$ then by \cite{S}  $P/JP$ has a free direct summand of rank one. Then by Nakayama's Lemma we can find $s\in J$ such that $P_{1+s}\cong Q\oplus B$, where $B=A_{1+s}$ is an affine algebra over $\ol{\mathbb{F}}_p$ and $Q$ is a projective $B$-module of rank $n-1$. Also note that $B[T]/IB[T]\cong (A[T]/I)_{1+s}$ gives us the fact $\dim(B[T]/IB[T])\le \dim(A[T]/I)$. Let $\alpha(T):P[T]\to I$ be any lift of $\ol{\phi}$, then we have a surjective map $\alpha(T)\otimes B[T]/ITB[T]=\ol{\phi}\otimes B[T]:(Q\oplus B)[T]\surj IB[T]/I^2TB[T]$. Let $\delta(T)$ be any lift of $\ol{\phi}\otimes B[T]$ and $f_0\in I$ be a monic polynomial. Then replacing $\delta(T)(0,1)$ by $\delta(T)(0,1)+T^kf_0^2$, for some suitably chosen $k>1$, we may assume that $\delta(T)(0,1)=f$ is a monic polynomial in $IB[T]$.
\paragraph{\textbf{Case - 1}} In this case we shall assume that $n\ge 3$. Define $C=B[T]/(J^2B[T],f)$. Then $C$ is an affine algebra over $\ol{\mathbb F}_p$ and since $(J^2B,f)$ contains a monic polynomial, namely, $f$, we have $\dim(C)=\dim(B/J^2B)=\dim(B/JB)= \dim(A/J)=\dim(A[T]/I)\le n-1$. $Q[T]\otimes C$ is a projective $C$-module of rank $n-1.$ Also note that $(\delta(T)\otimes C)|_{Q[T]\otimes C}=(\ol{\phi}\otimes C)|_{Q[T]\otimes C}:Q[T]\otimes C\surj IC/I^2TC $ is a surjective map with $\rank(Q[T]\otimes C)=n-1\ge \dim(C)$. Then by Lemma \ref{IMMTT} there exists a surjective map $\ol{\psi}:Q[T]\otimes C\surj IC $, which lifts $\ol{\phi}\otimes C$. Let $\psi:Q[T]\otimes B[T] \to IB[T]$ be a lift of $\ol{\psi}$. In the ring $B[T]$ we have $Im(\psi)+J^2B[T]+(f)=IB[T]$. Then by \cite{NMKL} there exists $e\in J^2B[T]$ with $e(1-e)\in (Im(\psi),f)$ such that $(Im(\psi),f,e)=I$. Define $I'=(Im(\psi),f,1-e)$, then $IB[T]\cap I'=(Im(\psi),f)$, $I'+J^2B[T]=B[T]$ and $I'$ contains a monic polynomial. Using \cite[Lemma 1.1, Chapter III]{Lam} we can find $t\in J^2B$ such that $1+t\in I'\cap B.$ Therefore, we get $IB_{1+t}[T]=(Im(\psi),f)B_{1+t}[T]$. Now note that since $t\in J^2B$, $t=\frac{b}{(1+s)^{k'}}$, for some $k'\ge 0$. Then by further localizing at $(1+s)^{k'}$  one can show that $(B_{1+t})=[(A_{1+s})_{1+t}]\cong A_{1+j}$, for some $j\in J$. Define $\omega'(T):P_{1+j}[T]\to I_{1+j}$, by $\omega'(T)|_{Q_{1+j}[T]}=\psi_{(1+s)}$ and $\omega'(T)(0,1)=f$. Then $\omega'(T)$ is a surjective map. Also note that  $\omega'(T) = {{(\delta}(T))_{1+j}}$ modulo $(J^2B_{1+j}[T],f,I^2TB_{1+j}[T])$. Since $f\in Im(\omega'_{1+s})$ is a monic polynomial, we can find a transvection $\zeta$ of $(Q_{1+j}[T]\oplus B_{1+j}[T])$ such that  $\omega'(T)|_{Q_{1+j}[T]}\zeta ={{(\delta}(T))_{1+j}}$ modulo $(J^2B_{1+j}[T],f,I^2TB_{1+j}[T])$. So we can replace $\omega'(T)$ by $\omega'(T)\zeta$ (without changing the notations) and may assume that $\omega'(T)\zeta = {(\delta(T))_{1+j}}=(\alpha(T))_{1+j}$ modulo $(J^2B_{1+j}[T],I^2TB_{1+j}[T])$. Define $\omega:P_{1+j}[T]\to I_{1+j}$, by $\omega(T)=\omega'(T)-(\omega'(0)-\alpha_{1+j}(0))$. Then we have the following:
\newline (1) $\omega(T)=\omega'(T)$ modulo $J^2B$.
\newline (2) $\omega(0)=\alpha_{1+j}(0).$ 

We get $\omega(T)=\alpha_{1+j}(T)=\ol{\phi}_{1=+j}$ modulo $I^2TB[T]$.
\par Again we have, $\alpha(0):P\to I(0)$ a surjective map. Therefore, $(\gamma(T)=)\alpha(0)\otimes A_{j}[T]:P_j[T]\to I(0)A_j[T](=A_j[T])$ is also a surjective map and $\gamma(T)=\ol{\phi} \mod(I^2TA_j[T])$.
\par Therefore, in the ring $A_{j(1+j)}[T]$ we have two surjective maps $(\omega(T))_j :P_{j(1+j)}[T]\surj I_{j(1+j)}$ and $(\gamma(T))_{1+j}:P_{j(1+j)}[T]\to I_{j(1+j)}$. Let $K_1=Ker(\omega(T))$ and $K_2=Ker(\gamma(T))$. Then we have the following two exact sequences
$$\begin{tikzcd}
0 \arrow{r} & (K_1)_j \arrow{r} & P_{j(1+j)}[T] \arrow{r}{\omega_j} & I_{j(1+j)} \arrow{r} & 0\\ 
0 \arrow{r} & (K_2)_{1+j} \arrow{r} & P_{j(1+j)}[T] \arrow{r}{\gamma_{1+j}} & I_{j(1+j)} \arrow{r} & 0.
\end{tikzcd}$$
\par Note that going modulo $T$ we have, $(\omega(0))_j=(\gamma(0))_{1+j}$, as they both matches with any lift of $\ol{\phi}_{j(1+j)}$ modulo $(I^2TA_{j(1+j)}[T])$. Since $I_{1+j}$ contains a monic polynomial then by \cite{Q} $K_1$ is extended. Also $K_2$ is extended from $A_j$ follows from the fact the map $\gamma$ itself is extended. Using \cite[Lemma 2]{P}  we can find an automorphism $\tau$ of $P_{j(1+j)}[T]$ such that $\tau(0)=Id$ and $\gamma_{1+j}\tau=\omega_j$. Then applying Quillen's Splitting Lemma \cite{Q} we get $\tau(T)=(\tau_1(T))_{1+j}(\tau_2(T))_{j}$, where $\tau_1(T)\in \Aut(P_{j}[T]) $ and  $\tau_2(T)\in \Aut(P_{1+j}[T]) $. Then a standard patching argument completes the proof in this case.

\paragraph{\textbf{Case - 2 }} In this case we shall assume that $n=2$. Using \cite{SM} we can assume that $\dim(A[T]/I)+1=2$. By Lemma \ref{NQS} it is enough to find $j\in J^2$ and a surjection $\omega:P_{1+j}[T]\surj I_{1+j}$ which lifts ${\ol{\phi}(T)}_{1+j}$.
\par  Since $I$ contains a monic polynomial, we have $\dim(A[T]/I)=\dim(A/J)=\dim(A/J^2)=1$. Let $C= \frac{A}{J^2}$,
and `bar' denote going modulo $J^2$, then in the ring $C[T]$ we have $$\ol \phi(T):\ol P[T]\surj \ol{I}/\ol{I^2T}.$$
Using Theorem \ref{IMMTT} we can find a lift $\ol{\psi}(T):P[T]\surj I/J^2[T]$ of $\ol{\phi}(T)$. Let us define $\omega(T)=\ol\psi(T)-\ol\psi(0)+\ol\phi(0):P[T]\surj I/J^2[T]$. Therefore, we get $I=\omega(P[T])+J^2[T]$.  Now since $I$ contains a monic polynomial and $J$ is a proper ideal in $A$, the ideal $\omega(P[T])$ must contains a monic polynomial. By \cite{NMKL} there exists $s\in J^2A[T]$ with $s(1-s)\in\omega(P[T])$. Let $I'=(\omega(P[T]),1-s>$. Then $I'$ contains a monic polynomial, $I\cap I'=\omega(P[T])$ and $I'+J^2[T]=A[T]$. Then by \cite[Lemma 1.1, Chapter III]{Lam} there exists $j\in J^2$ such that $1+j\in I'\cap A$. Therefore, we get $\omega(P[T])_{1+j}:P_{1+j}[T]\surj I_{1+j}$,  is a surjective lift of $\ol\phi(T)_{1+j}$ and this completes the proof. \qed

\bc
Let $A$ be an affine algebra over $\ol{\mathbb{F}}_p$ and $I\subset A[T]$ be any proper ideal containing a monic polynomial. Let $P$ be a projective $A$-module of rank $n$, where $n\ge \max\{ (\dim A[T]/I+1),2\}$ and $\lambda:P\surj I(0)$ be a surjection. Suppose that there exists a surjective map $\ol{\phi}:P[T]/IP[T]\surj I/I^2$ such that $\ol{\phi(0)}=\lambda(0)\otimes A/I(0)$. Then there exists a surjective map $\phi:P[T]\surj I$ such that $\phi$ lifts $\ol{\phi}$ and $\phi(0)=\lambda$.
\ec
\proof Follows from \cite[Remark 3.9 ]{BR} and using Theorem \ref{NQ}.\qed 

\bt \label{Mtrick}
Let $A$ be an affine algebra over  $\overline{\mathbb{F}}_p$. $I=I_1\cap I_2\subset A[T]$, is an ideal and $P$ is a projective $A$-module such that
\begin{enumerate}[\quad\quad(1)]
	\item $I_1$ contains a monic polynomial.
	\item $I_2 =I_2(0)A[T]$ is an extended ideal. 
	\item $I_1+I_2=A[T]$.
	\item $\rank(P)=n\ge \max\{(\dim(A[T]/I_1)+1),2\}$.
\end{enumerate}
Suppose that there exist surjections $\rho:P\surj I(0)$ and $\ol{\delta}:P[T]/I_1P[T]\surj I_1/I_1^2$
such that $\ol{\delta}=\rho\otimes A/I_1(0)$. Then there exists a surjection $\eta:P[T]\surj I$ such that $\eta(0)=\rho.$
\et
\proof Using Corollary \ref{QN}
and following the same argument used as in Theorem \ref{Mtrick1} the proof follows.\qed

\bibliographystyle{abbrv}


\begin{thebibliography}{10}
	
	\bibitem{af}
	A.~{A}sok and J.~{F}asel (with an appendix~by {M}. {K}.~{D}as).
	\newblock Euler class groups and motivic stable cohomotopy,.
	\newblock {\em {\it J. Eur. Math. Soc.} {\bf 24} (2022), 2775–2822}, 2022.
	
	\bibitem{bass}
	H.~{B}ass.
	\newblock Algebraic $K$-theory.
	\newblock W. A. Benjamin Inc. 1968.
	
	
	\bibitem{SMBmonic}
	S.~M. Bhatwadekar.
	\newblock Inversion of monic polynomials and existence of unimodular elements.
	{II}.
	\newblock {\em Mathematische Zeitschrift}, 200(2):233--238, 1989.
	
	\bibitem{MKDSMB}
	S.~M. Bhatwadekar and M.~K. Das.
	\newblock Projective generation of curves ({III}).
	\newblock {\em International Mathematics Research Notices}, 2015(4):960--980,
	oct 2013.
	
	\bibitem{BHMK}
	S.~M. Bhatwadekar and M.~K. Keshari.
	\newblock A question of {N}ori: {P}rojective generation of ideals.
	\newblock {\em K-Theory}, 28(4):329--351, apr 2003.
	
	\bibitem{SMBHLRR}
	S.~M. Bhatwadekar, H.~Lindel, and R.~A. Rao.
	\newblock The {B}ass-{M}urthy question: Serre dimension of {L}aurent polynomial
	extensions.
	\newblock {\em Inventiones Mathematicae}, 81(1):189--203, feb 1985.
	
	\bibitem{BR}
	S.~M. Bhatwadekar and R.~Sridharan.
	\newblock Projective generation of curves in polynomial extensions of an affine
	domain and a question of {N}ori.
	\newblock {\em Inventiones Mathematicae}, 133(1):161--192, jun 1998.
	
	\bibitem{BRS}
	S.~M. Bhatwadekar and R.~Sridharan.
	\newblock Zero cycles and the {E}uler class groups of smooth real affine
	varieties.
	\newblock {\em Inventiones Mathematicae}, 136(2):287--322, apr 1999.
	
	\bibitem{SMBB3}
	S.~M. Bhatwadekar and R.~Sridharan.
	\newblock The {E}uler class group of a {N}oetherian ring.
	\newblock {\em Compositio Mathematica}, 122(2):183--222, 2000.
	
	\bibitem{BRS01}
	S.~M. Bhatwadekar and R.~Sridharan.
	\newblock {O}n a question of {R}oitman.
	\newblock {\em J. Ramanujan Math. Soc}, 16(1):45--61, 2001.
	
	\bibitem{BR3}
	S.~M. Bhatwadekar and R.~Sridharan.
	\newblock On {E}uler classes and stably free projective modules.
	\newblock In {\em Algebra, arithmetic and geometry, {P}art {I}, {II} ({M}umbai,
		2000)}, volume~16 of {\em Tata Inst. Fund. Res. Stud. Math.}, pages 139--158.
	Tata Inst. Fund. Res., Bombay, 2002.
	
	\bibitem{MKD1}
	M.~K. Das.
	\newblock The {E}uler class group of a polynomial algebra.
	\newblock {\em Journal of Algebra}, 264(2):582--612, jun 2003.
	
	\bibitem{d2}
	M.~K. Das.
	\newblock Revisiting {N}oris question and homotopy invariance of {E}uler class
	groups.
	\newblock {\em Journal of K-Theory 8 (2011), 451-480}, 2011.
	
	\bibitem{Das2012}
	M.~K. Das.
	\newblock On triviality of the {E}uler class group of a deleted neighbourhood
	of a smooth local scheme.
	\newblock {\em Transactions of the American Mathematical Society},
	365(7):3397--3411, dec 2012.
	
	\bibitem{MKD}
	M.~K. Das.
	\newblock On a conjecture of {M}urthy.
	\newblock {\em Advances in Mathematics}, 331:326--338, 2018.
	
	\bibitem{MKDIMRN}
	M.~K. Das.
	\newblock A criterion for splitting of a projective module in terms of its
	generic sections.
	\newblock {\em International Mathematics Research Notices},
	2021(13):10073--10099, jun 2019.
	
	\bibitem{DTZ}
	M.~K. Das, S.~Tikader, and M.~A. Zinna.
	\newblock Orbit spaces of unimodular rows over smooth real affine algebras.
	\newblock {\em Inventiones mathematicae}, 212(1):133--159, oct 2017.
	
	\bibitem{DTZ2}
	M.~K. Das, S.~Tikader, and M.~A. Zinna.
	\newblock "{${ {P}}^1$}-gluing'' for local complete intersections.
	\newblock {\em Mathematische Zeitschrift}, 294(1-2):667--685, apr 2019.
	
	\bibitem{EE}
	D.~Eisenbud and E.~G. Evans, Jr.
	\newblock Generating modules efficiently: theorems from algebraic {$K$}-theory.
	\newblock {\em Journal of Algebra}, 27:278--305, 1973.
	
	\bibitem{MPHIL}
	M.~K. Keshari.
	\newblock Euler class group of a {N}oetherian ring.
	\newblock {\em M-Phil Thesis-2001 (52 page).}, Aug. 2001.
	
	\bibitem{STMKK}
	M.~K. Keshari and S.~Tikader.
	\newblock On a question of {M}oshe {R}oitman and {E}uler class of stably free
	module.
	\newblock Feb. 2022 (available at https://arxiv.org/abs/2202.06291).
	
	\bibitem{AS}
	A.~Krishna and V.~Srinivas.
	\newblock Zero cycles on singular varieties.
	\newblock In {\em Algebraic cycles and motives. {V}ol. 1}, volume 343 of {\em
		London Math. Soc. Lecture Note Ser.}, pages 264--277. Cambridge Univ. Press,
	Cambridge, 2007.
	
	\bibitem{NMK}
	N.~M. Kumar.
	\newblock On two conjectures about polynomial rings.
	\newblock {\em Inventiones Mathematicae}, 46(3):225--236, oct 1978.
	
	\bibitem{Lam}
	T.~Y. Lam.
	\newblock {\em Serre's problem on projective modules}.
	\newblock Springer Monographs in Mathematics. Springer-Verlag, Berlin, 2006.
	
	\bibitem{MM}
	S.~Mandal and M.~Pavaman~Murthy.
	\newblock Ideals as sections of projective modules.
	\newblock {\em Journal of the Ramanujan Mathematical Society}, 13(1):51--62,
	1998.
	
	\bibitem{MR}
	S.~Mandal and R.~Sridharan.
	\newblock Euler classes and complete intersections.
	\newblock {\em Journal of Mathematics of Kyoto University}, 36(3):453--470,
	1996.
	
	\bibitem{mv}
	S.~Mandal and P.~L.~N. Varma.
	\newblock On a question of {N}ori: the local case.
	\newblock {\em Communications in Algebra}, 25(2):451--457, jan 1997.
	
	\bibitem{SM}
	S.~{M}andal (with an appendix by~{M}adhav {V}.~{N}ori).
	\newblock Homotopy of sections of projective modules.
	\newblock {\em J. Algebraic Geom. 1 (1992)}, (4):639--646, 1992.
	
	\bibitem{NMKL}
	N.~Mohan~Kumar.
	\newblock Complete intersections.
	\newblock {\em Journal of Mathematics of Kyoto University}, 17(3):533--538,
	1977.
	
	\bibitem{AMM}
	N.~Mohan~Kumar, M.~P. Murthy, and A.~Roy.
	\newblock A cancellation theorem for projective modules over finitely generated
	rings.
	\newblock In {\em Algebraic geometry and commutative algebra, {V}ol. {I}},
	pages 281--287. Kinokuniya, Tokyo, 1988.
	
	\bibitem{P}
	B.~Plumstead.
	\newblock The conjectures of {E}isenbud and {E}vans.
	\newblock {\em American Journal of Mathematics}, 105(6):1417--1433, 1983.
	
	\bibitem{pop}
	D.~Popescu.
	\newblock Polynomial rings and their projective modules.
	\newblock {\em Nagoya Mathematical Journal}, 113:121--128, mar 1989.
	
	\bibitem{Q}
	D.~Quillen.
	\newblock Projective modules over polynomial rings.
	\newblock {\em Inventiones Mathematicae}, 36:167--171, 1976.
	
	\bibitem{S}
	J.-P. Serre.
	\newblock Modules projectifs et espaces fibrés à fibre vectorielle.
	\newblock {\em Séminaire Dubreil. Algèbre et théorie des nombres},
	11(2):1--18, 1957-1958.
	
	\bibitem{R}
	R.~G. Swan.
	\newblock Serre's problem.
	\newblock pages 1--60. Queen's Papers on Pure and Applied Math., No. 42.
	Conference on {C}ommutative {A}lgebra -- 1975 ({Q}ueen's {U}niv., {K}ingston,
	{O}nt., 1975), 1975.
	
	\bibitem{VdKM}
	W.~van~der Kallen.
	\newblock A module structure on certain orbit sets of unimodular rows.
	\newblock {\em Journal of Pure and Applied Algebra}, 57(3):281--316, apr 1989.
	
	\bibitem{MNLvdK}
	W.~van~der Kallen.
	\newblock From {M}ennicke symbols to {E}uler class groups.
	\newblock {\em Proceedings of the International Colloquium on Algebra,
		Arithmetic and Geometry. Mumbai 2000, Part II, 341-354}, Oct. 2000.
	
	\bibitem{V}
	L.~N. Vaser{\v{s}}te{\u{\i}}n.
	\newblock Stabilization of unitary and orthogonal groups over a ring with
	involution.
	\newblock {\em Mathematics of the {USSR}-Sbornik}, 10(3):307--326, apr 1970.
	
	\bibitem{SV}
	L.~N. Vaser{\v{s}}te{\u{\i}}n and A.~A. Suslin.
	\newblock {S}erre's problem on projective modules over polynomial rings, and
	algebraic {K}-theory.
	\newblock {\em Mathematics of the {USSR}-Izvestiya}, 10(5):937--1001, oct 1976.
	
\end{thebibliography}

\end{document}